\documentclass[a4paper, a4wide,11pt]{amsart} 
\usepackage[utf8]{inputenc}
\usepackage{esint} 
\usepackage{amsmath,amsthm, amssymb}
\usepackage{color}
\usepackage{graphicx}
\usepackage{geometry}
\usepackage{mathrsfs}
\newtheorem{theorem}{Theorem}[section]
\newtheorem{lemma}[theorem]{Lemma}
\newtheorem{prop}[theorem]{Proposition}

\newtheorem{defn}[theorem]{Definition}
\newtheorem{cor}[theorem]{Corollary}
\newtheorem{rem}[theorem]{Remark}

\usepackage{enumerate}
\usepackage{abraces} 

\usepackage[T3,T1]{fontenc}
\DeclareSymbolFont{tipa}{T3}{cmr}{m}{n}
\DeclareMathAccent{\invbreve}{\mathalpha}{tipa}{16}

\makeatletter
\@namedef{subjclassname@2020}{\textup{2020} Mathematics Subject Classification}
\makeatother

\makeatletter
\def\widebreve{\mathpalette\wide@breve}
\def\wide@breve#1#2{\sbox\z@{$#1#2$}%
     \mathop{\vbox{\m@th\ialign{##\crcr
\kern0.08em\brevefill#1{0.8\wd\z@}\crcr\noalign{\nointerlineskip}%
                    $\hss#1#2\hss$\crcr}}}\limits}
\def\brevefill#1#2{$\m@th\sbox\tw@{$#1($}%
  \hss\resizebox{#2}{\wd\tw@}{\rotatebox[origin=c]{90}{\upshape(}}\hss$}
\makeatletter

\DeclareMathOperator*\Id{Id}
\DeclareMathOperator*\dist{dist}

\DeclareMathOperator*\Om{\Omega}
\DeclareMathOperator*\diam{diam}

\usepackage[colorinlistoftodos,prependcaption,textsize=tiny]{todonotes}

\newcommand*\R{\mathbb{R}}
\newcommand*\tr{{\rm Tr}}
\newcommand\ds{\displaystyle}

\newcommand\numberthis{\addtocounter{equation}{1}\tag{\theequation}}

\numberwithin{equation}{section}
\numberwithin{figure}{section}

\title{The biharmonic optimal support problem}
 
  \author[A. Lemenant]{Antoine Lemenant}
 \author[M. R. Pakzad]{Mohammad Reza Pakzad}

\address[M. R. Pakzad]{Laboratoire Imath, Université de Toulon, CS 60584 - 83041, TOULON CEDEX 9, France}
\email{\vspace{-0.15in}pakzad@univ-tln.fr}
  
 \address[A. Lemenant]{Institut universitaire de France (IUF) et Universit\'e de Lorraine, CNRS, IECL, F-54000 Nancy, France}
\email{antoine.lemenant@univ-lorraine.fr}


\date{}

\begin{document}

 \begin{abstract} { We establish  a $\Gamma$-convergence result for  $h\to 0$  of  a thin nonlinearly elastic 3D-plate of thickness $h>0$ which  is assumed to be  glued to a support region in the 2D-plane  $x_3=0$ over the $h$-2D-neighborhood of a given closed set $K$.  In the regime of very small  vertical forces we identify the $\Gamma$-limit as being the bi-harmonic energy, with Dirichlet condition on the gluing region $K$, following a general strategy by Friesecke,  James, and M\"uller that we have to adapt in presence of the glued region. Then we introduce a shape optimization problem that we call ``optimal support problem'' and which aims to find the best glued plate. In this problem  the bi-harmonic energy is optimized among all possible glued regions $K$ that we assume to be connected and for which we penalize the length. By relating the dual problem with  Griffith almost-minimizers, we are able to prove that any minimizer   is $C^{1,\alpha}$ regular outside a set of Hausdorff dimension strictly less then one.   } 
  \end{abstract} 
  
  \maketitle
  \vspace{1cm}
\tableofcontents


\newpage
\section{Introduction}

Let $\Omega \subset \mathbb{R}^2$ a bounded Lipschitz domain, and  $K\subset \overline{\Omega}$ be a closed set.  For any $k\geq 1$ we consider the subspace of $W^{k,2}(\Omega)$ defined by 
$$H^k_{0,K}(\Omega):=\left\{u \in W^{k,2}(\Omega) \text{ such that } \partial^{\alpha}u|_{K}=0  \text{ for all multiindices $\alpha$ such that } 0\leq \alpha \leq k-1\right\},$$
where  by $\partial^{\alpha}u|_K =0$ we mean that $ \lim_{r\to 0} \left( \fint_{B_r(x_0) \cap \Omega}|\partial^\alpha u| \; dx\right) = 0 \; \text{ for  $(k-|\alpha|,2)$-q.e. } x_0 \text{ on }K.$

It follows from the literature (see Lemma \ref{QeLemma0} in the appendix for a proof), that if $\Omega \subset \R^N$ is a Lipschitz domain then
 $$H^k_{0,K}(\Omega)=\Big\{u \in W^{k,2}(\Omega) \text{ such that } \mathcal{E}(u) \in W^{k,2}_0(\R^N \setminus K)\Big\},$$
 where $\mathcal{E}$ is any extension operator  $\mathcal{E}:W^{k,2}(\Omega)\to W^{k,2}(\R^N)$ and that $H^k_{0,K}(\Omega)$ stands for a closed subspace of $W^{k,2}(\Omega)$.

Letting $f\in L^2(\Omega)$, this provides the existence of  $u_K$ being the unique solution of 
\begin{eqnarray}\label{complianceprob}
\left\{
\begin{array}{c}
\Delta^2 u_K = f\\
u\in {H^2_{0,K}(\Omega).}
\end{array}
\right.
\end{eqnarray}

In particular, $u_K$ is the unique  critical point to the  bi-harmonic energy

\begin{equation}\label{lin-comp}
  E(u):=\frac 12 \int_\Omega |\nabla ^2 u |^2  \; dx - \int_{\Om} uf \; dx  
\end{equation} over {$H^2_{0,K}(\Omega)$},
 and the associated  compliance  energy  is
\begin{eqnarray}
\int_{\Omega} u_K f \; dx=\int_{\Omega} |\nabla^2 u_K|^2  \;dx. \label{Compliance2}
 \end{eqnarray}

In this paper we intend to derive the functional \eqref{lin-comp} as the $\Gamma$-limit of a support-gluing problem for the 3d nonlinear thin elastic body subject to vertical body forces very small   with respect to its thickness, which is the  main result of the present paper.

To be more precise, we denote  by $\Om^h:= \Om \times (0,h)$, which models a thin elastic  plate which  is assumed to be  glued to a support
in the plane $x_3=0$ by applying a surface glue to the support region
$$
K_h:= \{ (x',0) \in  \Omega \times \{0\};\,\, \dist (x', K) < h\}.
$$ 
Then we assume the plate to be  subject to a  vertical body force ${\bf f}^h: \Om^h \to \R^3$, ${\bf f}^h:= (0,0, h^\alpha \tilde f)$, where $\alpha > 2$,  
 $\tilde f(x',x_3) = f(x')$ for   $f\in L^2(\Omega)$.

In the most simplified framework (see Section \ref{section1} for a general case), the bulk elastic energy of a deformation 
 $y: \Om^h \to \R^3$ is given by 
 $$
  E^h(y)   = \frac 1h  \int_{\Omega^h} \frac{1}{2}\dist^2(\nabla y, SO(3))\;dx.
 $$  
  The deformation of the plate $\Om^h$ subject to body forces and gluing constraint is then variationaly modeled by minimizing the energy functional 
 $$
J^h(y):= E^h(y)  -  \frac 1h \int_{\Om^h} {\bf f}^h \cdot y \; dx,
 $$ among all deformations in the class
\begin{equation}\label{glued-intro}
 \mathcal {A}^{h}_K := \{y\in W^{1,2}(\Omega^h, \R^3);\,\, y|_{K_h \times \{0\}}  = {\rm id}_{K_h \times \{0\}}\}.
\end{equation} 
The first goal of this paper is to establish the $\Gamma$-convergence of $J^h$ towards the bi-harmonic functional \eqref{lin-comp}, in the regime of  very small deformations in the sense that there exists $\beta>4$ such that 
 \begin{eqnarray}
\displaystyle \limsup_{h \to 0}  \frac 1{h^\beta} \Big  (J^h (y^h) - \inf_{\mathcal {A}^{h}_K} J^h \Big ) =0. \label{boundedE}
\end{eqnarray}

Here is { a shorter version} of one of the main statements of our paper, {stated, { moreover,} for simplicity here in the introduction, for the most simplified functionals $E^h$ and $J^h$ as above}. In order to state this result we   introduce   the average vertical displacement   functions
 $$
 u_h\in  W^{1,2} (\Omega,\R), \quad  U_h:= \frac 1h \int_0^h (y^h_3 (\cdot, x_3) - hx_3)   \, dx_3,\quad  u_h:= \frac{1}{h^{\beta/2-1}} U_h,
 $$ {where $y^h_3$ is the vertical component of the deformation $y^h$}. 
%

\begin{theorem}[$\Gamma$-convergence]\label{main11}Assume that $\Omega \subset \R^2$ is a Lipschitz domain and $K\subset \overline{\Omega}$ is a closed set satisfying ${\rm Cap}_{1,2}(K)>0$. Assume furthermore that  $\beta>4$. Then the functional $J^h$ $\Gamma$-converges to the energy in \eqref{lin-comp} in the following sense. 

\begin{enumerate}
\item for any sequence $h\to 0$ such that $y^h \in \mathcal{A}^h_K$ with energy bound $E^h(y^h) \precsim h^\beta$, there eone can find {$u\in H^2_{0,K}(\Omega)$  and} a  subsequence such that   $\ds u_h  \xrightarrow{\mbox{\tiny in} \,\, W^{1,2}} u \in W^{2,2}(\Omega)$,  and
$$\liminf_h \frac{1}{h^\beta} J^h(y^h) \geq  \int_{\Omega} \frac{1}{24}|\nabla^2 u |^2 - u f \; dx.$$

\item for every $u\in H^2_{0,K}(\Omega)$, there exists a sequence $y^h \in  \mathcal{A}^h_K$  such that $\ds u_h  \xrightarrow{\mbox{\tiny in} \,\, W^{1,2}} u $ and 
$$\limsup_h \frac{1}{h^\beta} J^h(y^h) \leq  \int_{\Omega} \frac{1}{24}|\nabla^2 u |^2 - u f \; dx.$$
   \end{enumerate}
\end{theorem}

The first part of Theorem \ref{main11} follows from Theorem \ref{compactness} which gives a more general compactness result for sequences satisfying  $E^h(y^h) \precsim h^\beta$. In Section \ref{liminf} we explain how to obtain the liminf inequality in (1)  (see Theorem \ref{liminf-thm}). The second part is the purpose of  Section \ref{recovery}. Notice that our main results are actually more general and work for a nonlinear energy of of the form $\int W(\nabla y) \;dx$ with $W$ satisfying some standard assumptions.

 For the bulk of the argument, the  proof or our result is inspired by  the seminal papers by Friesecke,  James, and M\"uller \cite{FJMrigid, FJMhier}, specially regarding the so called linearized von K\'arm\'an theory in \cite[Theorem 2]{FJMhier}. The novelty we need to take care of here is the \lq \lq gluing'' part  of the problem  on the support $K_h$ and its passage to the limit. Notice that we only assume $y={\rm id}$ on the 2D bottom part $K_h \subset \Omega \times \{0\}$ which differs from a standard Dirichlet or clamped boundary conditions. In particular, a key step to obtain compactness of a sequence with bounded energy in a thin plate is to be able to approximate it by piecewise constant maps with values in $SO(3)$. In our context, we have to adapt this approximation by insuring the map to be constant equal to identity on the support set $K_h$. As $K_h$ is a very thin $2D$ set in $\R^3$, this part is non trivial and uses thin properties on Sobolev functions  (see Corollary \ref{for-R-cor} for a statement).

 Since our goal is to arrive to the bi-harmonic problem,  we consider in this paper only the case of $\beta>4$, even if some other regimes could also be probably investigated.  As a result, we can interpret the solution $u_K$  in \eqref{complianceprob} as modeling in the linear regime a vertical displacement of the 2D-plate $\Omega$,  attached, or supported, onto the set $K$, and by \eqref{Compliance2} the energy associated to this displacement.  
 Then,  one can  seek for the ``best way'' of attaching the plate, when the support $K$ is penalized by its length. This leads to the following bi-harmonic optimal support problem:
\begin{eqnarray}
\min_{K \subset \overline{\Omega}} \quad  \int_{\Omega} |\nabla^2 u_K   |^2 \; dx + \mathcal{H}^1(K), \label{probmain2}
\end{eqnarray} 
where the minimum is taken over all compact and connected subsets $K$.\\

It is worth mentioning that replacing the bi-harmonic operator by the standard Laplace operator would lead to the so-called   ``optimal compliance problem'' that was studied before in many papers (see \cite{Approx,Butazzo-Santambrogio,Opt,MR3063566,MR3195349,B1,B2,B3}).   A variant with a $p$-Laplacian has been also studied in \cite{B1,B2,B3}. 

In  other words, our paper deals with a bi-laplacian variant of the standard optimal compliance problem, which turns out to be relevant from the mechanical point of view. One of the motivation for this paper was indeed to provide a better justification of the so-called ``optimal compliance problem''. In the second part of the paper we prove existence and regularity results for a minimizer $K$ of problem \eqref{probmain2}, leading to a  second   result proved within this paper.

 \begin{theorem} \label{main22} For every bounded domain $\Omega \subset \R^2$  there exists a minimizer $K$ for the optimal support problem in \eqref{probmain2}. Moreover, if $K\cup \partial \Omega$ is connected, then any minimizer $K$ is an almost minimizer of the Griffith functional, and therefore is locally $C^{1,\alpha}$ regular inside $\Omega$, except for a  singular set of points with Hausdorff dimension strictly less then $1$.
 \end{theorem}

 To prove Theorem \ref{main22} we consider the dual formulation of the problem in \eqref{probmain2}, and prove that the minmiizers of the dual problem are almost minimizers of the so-called Griffith functional, already studied in the literature. By use the regularity results contained in  \cite{lem-lab} and \cite{lem-lab0}, we obtain the conclusion. 
 
\subsection {Acknowledgments.}
This project was based upon work supported by the National Science Foundation. M.R.P. was supported by the NSF award DMS-1813738.  This project was completed while the author was visiting Institut \'Elie Cartan de Lorraine at Nancy supported by CNRS funding, and he is grateful for their respective hospitality and support.  A. L. is partially supported by the project ANR-18-CE40-0013 SHAPO financed by the French Agence Nationale de la Recherche (ANR).


   \section{Justification of  the  biharmonic support model from 3d nonlinear elasticity}
 
 \subsection{The 3d support-gluing model}
 \label{section1}

As before in the introduction, let $\Omega \subset \mathbb{R}^2$ be a bounded Lipschitz domain, $K\subset \overline \Omega$ be a closed set and let $\Om^h:= \Om \times (0,h)$ a thin elastic  plate. As before it is assumed that the thin plate is glued to a support
in the plane $x_3=0$ by applying a surface glue to the support region
$$
K_h:= \{ (x',0) \in  \Omega \times \{0\};\,\, \dist (x', K) < h\},
$$ 

 and is subject to a  vertical body force ${\bf f}^h: \Om^h \to \R^3$, ${\bf f}^h:= (0,0, h^\alpha \tilde f)$, where $\alpha > 2$,  
 $\tilde f(x',x_3) = f(x')$ for   $f\in L^2(\Omega)$. 
 Here we have opted for a simplified set of assumptions on the body forces in order to focus on the new contribution.   
   
In what follows  
$$
SO(n):= \{ R\in  \R^{n\times n}; \,\, R^T R = \Id, \,\, \det R>0 \}    
 $$ is the special orthogonal group of  3d rotation matrices.   The elastic density or potential $W: \R^{3\times 3}\to \R $ 
 is assumed to satisfy the following natural conditions for all 
 $F\in \R^{3\times 3}$:

 \begin{itemize}
 \item Normalization: $W(F)\ge 0, \,\, W(\Id)=0$.
  \item Frame invariance: $\forall R \in SO(3), W(RF) = W(F)$.
 \item Non-degeneracy: $W(F) \ge c \dist^2 (F, SO(3))$ for a constant $c>0$.
  \end{itemize} Moreover we will assume that $W$ is of class $C^2$ in a neighborhood of $SO(3)$. The bulk elastic energy of a deformation 
 $y: \Om^h \to \R^3$ is given by 
 $$
  E^h(y)   = \frac 1h  \int_{\Omega^h} W(\nabla y) \;dx.
 $$  The deformation of the plate $\Om^h$ subject to body forces and gluing constraint is then variationally modeled by minimizing the energy functional 
 $$
J^h(y):= E^h(y)  -  \frac 1h \int_{\Om^h} {\bf f}^h \cdot y \; dx,
 $$ among all deformations in the class
\begin{equation}\label{glued}
 \mathcal {A}^{h}_K := \{y\in W^{1,2}(\Omega^h, \R^3);\,\, y|_{K_h \times \{0\}}  = {\rm id}_{K_h \times \{0\}}\}.
\end{equation} 

 Notice that we will still denote by $E^h$ and $J^h$ some functionals that was introduced before in the introduction in the particular  case of  $W(F)=\frac{1}{2} \dist^2 (F, SO(3))$.
 
 
 For $\beta >0$,  we say the sequence $y^h$ is a $\beta$-minimizing sequence for $J^h$ whenever
 $$
\displaystyle \limsup_{h \to 0}  \frac 1{h^\beta} \Big  (J^h (y^h) - \inf_{\mathcal {A}^{h}_K} J^h \Big ) =0.
 $$
  
   \subsection{The limiting 2D energy}

 Following \cite{FJMhier} we introduce the linearized energy of second order  
 $$Q_3(F):=\frac{\partial^2W}{\partial F^2}(Id)(F,F),$$
 and $Q_2: \mathbb{R}^{2\times 2}\to \R$ defined by 
 $$Q_2(G)= \min_{a \in\R^3} Q_3(G+a\otimes e_3 +e_3\otimes a).$$
 Under our assumptions, both forms are positive, semidefinite, convex and positive definite on symmetric matrices.
 
 For the special case of isotropic elasticity, i.e. when $W(FR) = W(F)$ for all $R\in SO(3)$ and $F\in \R^{3\times 3}$, it can be shown that 
   $$Q_3(F)=2\mu \left| \frac{F+F^T}{2} \right|^2  +\lambda (\tr F)^2,$$
 and  
  $$Q_2(G)= 2\mu \left| \frac{G+G^T}{2} \right|^2  +\frac{2\mu \lambda}{2\mu+\lambda} (\tr G)^2,$$ where $\mu >0$ and $\lambda\geq 0$. In particular if $W(F) = \frac 12 {\rm dist}^2 (F, SO(3))$, then $\mu=1/2$ and  $\lambda=0$.

   \subsection{Friesecke-James-M\"uller rigidity estimate in presence of affine boundary conditions}
   In this section we present a corollary of the celebrated geometric rigidity estimate of Friesecke-James-M\"uller
  \cite{FJMrigid} which will be a necessary ingredient of the compactness argument. We first state this rigidity estimate:
  
  \begin{theorem}\cite[Theorem 3.1]{FJMrigid}\label{rig-es-th}
Let $n\ge 2$ and $\Omega\subset \R^n$ be a bounded Lipschitz domain.  Then there exists a constant $C= C(\Omega)$,  such that for all mapping  $u\in W^{1,2}(\Omega, \R^n)$, there is a rotation $R\in SO(n)$ with 
\begin{equation}\label{rig-es}
 \|\nabla y - R\|^2_{L^2 (\Omega)} \le C   \|{\dist} (\nabla y , SO(n))\|^2_{L^2(\Omega)} 
  \end{equation}
  \end{theorem}  
  
 We would like to identify the \lq \lq average rotation" $R$  in \eqref{rig-es} based on the boundary conditions  applied to $y$ 
 on portions of $\partial \Omega$. In particular, we expect that the rotation $R$ can be chosen to be uniformly  the identity matrix 
 $\Id$ for all  mappings for which $y|_S = {\rm id}_S$ for an open set $S\subset \partial \Omega$: 
 
 \begin{defn}
 We say that a matrix $F\in \R^{n\times n}$ is rank-one connected to $SO(n)$ whenever there exists $R\in SO(n)$ for which 
 ${\rm rank} (F-R)=1$. Note that no element of $SO(n)$ itself  does enjoy this property.
 \end{defn}
 
\begin{cor}\label{fixed-Id}
Let $n\ge 2$ and $\Omega\subset \R^n$ be a bounded Lipschitz domain. For $F\in \R^{n\times n}$, and $S$ an open connected subset of $\partial  \Omega$, we define
$$
\mathcal{A}_{S, F}:= \{y\in W^{1,2}(\Omega, \R^n);\,\, y|_S (x)= Fx\}. 
$$   Let  $R_F$ be any closest element of $SO(n)$ to $F$.

(i) If $F$ is not rank-one connected to $SO(n)$,  there exists $C=C(\Omega, S, F)$ such that for all $y\in \mathcal{A}_{S, F}$, 
 $$
 \|\nabla y - R_F\|^2_{L^2 (\Omega)} \le C  \|{\dist} (\nabla y , SO(n))\|^2_{L^2(\Omega)}. 
 $$ 

(ii) Assume  moreover  that $S$ is inside no hyperplane of $\R^n$. Then, there exists $C=C(\Omega, S)$ such that for all $y\in \mathcal{A}_{S, F}$, 
 $$
 \|\nabla y - R_F\|^2_{L^2 (\Omega)} \le C \|{\dist} (\nabla y , SO(n))\|^2_{L^2(\Omega)}. 
  $$ 
 \end{cor} 

\begin{rem} The estimate can fail to be true if  ${\rm rank}(F-R)=1$ and $S$ is a subset of a hyperplane.  
As a counter-example, take such $R, F$ in a manner that $R \neq R_F$, let $y=Rx$, and choose  $S$ such that $Rx=Fx$ on S.
 \end{rem}

  \begin{proof} The proof follows  the same approach as in \cite[Proposition 3.4]{DNP2002}. We write the full details for the reader's convenience.  In what follows the constant $C$ might differ from line to line but it will always depend on $\Omega$. Its dependance on $S$ or $F$ will be clarified on the way.  
  
  Consider an arbitrary point $p\in S \subset \R^n$ to be fixed later. Since the domain is Lipschitz, there is a hyperplane $P\ni p$ such that  $\partial \Omega$ is locally the graph of a Lipschitz function over an open subset $U\subset P$ containing $p$. Let $V$ be the intersection of the projection of $S$ on $P$ with $U$. Fix a rotation $R_0 \in SO(n)$ such that $\psi(x) := R_0 x + p$  maps $\R^{n-1}\times \{0\}$  to $P$.   We let $\widetilde \Omega:= \psi ^{-1} (\Omega)$, $\tilde S := \psi^{-1}(S) \subset \partial \widetilde \Omega$ and we fix $r>0$ for which $\psi (B_r(0) \times \{0\})\subset V$, where $B_r(0)$ is the ball of radius $r$ in $\R^{n-1}$. Hence there exists a  Lipschitz function $g: B_r(0) \to \R$ such that  whose graph is an open subset of $\tilde S$.  We define $\phi_{}: B_r(0) \to \tilde S$ to be the graph parameterization of the portion of $\tilde S$ which is over $B_r(0)$:
 $$
 \quad   \phi_{}(z) := \left [ \begin{array}{c} z \\ g( z)    \end{array}  \right ]  = 
\left [ \begin{array}{c} z \\ \tilde g( z)    \end{array}  \right ] +   \left [  \begin{array}{c} 0 \\  c    \end{array}  \right ],  
 $$ where 
 $$
\ds  \tilde g(z):= g(z) - \fint_{B_r(0)} g(z) dz , \quad c:=   \fint_{B_r(0)} g(z)\, dz .     
 $$  
  
  \medskip
In the context of part (ii) we can choose $p$ and $r$ such that $\tilde g$ is not a linear function over $B_r(0)$, for otherwise $S$, being connected, would become part of a hyperplane.
   
   \medskip
Now assume that $y\in \mathcal{A}_{S, F}$ is given and let
$$
E:= \|{\dist} (\nabla y , SO(n))\|^2_{L^2(\Omega)}. 
$$ Applying Theorem \ref{rig-es-th} to $y$ on $\Om$ we obtain that for a uniform constant $C>0$, there exists 
$R\in SO(n)$ for which 
\begin{equation}\label{rig-es-orig}
\|\nabla   y - R\|^2_{L^2(\Om)} \le  C E,
\end{equation} We claim that under assumptions of parts (i) or (ii), $y\in \mathcal{A}_{S, F}$ implies 
\begin{equation}\label{fin-es}
|F-R|^2\le CE,
\end{equation} with $C$ depending on the respective claimed variables. Let us observe that this is sufficient to conclude the proof of the corollary: Indeed, our claim yields in view of the definition of $R_F$ and the fact that $R \in SO(n)$:
$$
|R_F - R | \le |F-R_F| + |F- R| \le 2 |F- R|  \le C\sqrt E.  
$$ This, combined with \eqref{rig-es-orig} implies
$$
\|\nabla  y - R_F\|^2_{L^2(\Om)} \le 2 ( \|\nabla y - R\|^2_{L^2(\Om)} +  |R_F - R|^2 \mathcal{L}^n(\Omega)) \le CE,
$$ which is the desired estimate. Hence what remains is to prove \eqref{fin-es}.

Letting $y:= y\circ \psi$ on $\widetilde \Omega$, and in view of $\nabla \tilde y = (\nabla y \circ \psi) R_0$, we obtain  from \eqref{rig-es-orig}
\begin{equation}\label{rig-es-tilde}
\|\nabla  \tilde y - RR_0\|^2_{L^2(\widetilde \Om)} \le  C E,
\end{equation}
We let $\ds b:= \fint_{\Om} (\tilde y (x)- RR_0 x)\; dx$. Applying the Poincar\'e inequality on $\widetilde \Omega$,  and the trace embedding of $W^{1,2} (\widetilde \Omega)$ into $L^2(\partial \widetilde \Omega)$, we obtain
$$
\|\tilde y - (RR_0 x+b) \|^2_{L^2(\tilde S)} \le C \|\tilde y- (RR_0 x +b)\|^2_{W^{1,2}(\widetilde \Omega)}\le C E. 
$$   But $\tilde y|_{\tilde S} = F(R_0x + p) = FR_0 x + Fp$.  This yields 
\begin{equation}\label{til-S-es}
 \|(F - R)R_0x + (Fp -b)\|^2_{L^2(\tilde S)} \le CE. 
\end{equation}  

 In what follows, and for any $F\in \R^{n\times n}$, we will denote by $\widehat F \in \R^{n \times (n-1)}$ and $F^{(n)} \in \R^n$, respectively the matrix made by the first $n-1$ columns of $F$, and the last column of $F$. \eqref{til-S-es} gives
  \begin{equation}\label{e-0}
 e:=  \|(F-R)R_0x + (Fp -b)\|^2_{L^2(\phi_{}(B_r(0))} \le  \|(F-R)R_0x + (Fp -b)\|^2_{L^2(\tilde S)} \le CE. 
\end{equation} We estimate from below the left hand side through change of variable:
$$
\begin{aligned}
e
& = \int_{\phi_{}(B_r(0))} \Big |(F-R)R_0x + (Fp -b)\Big |^2 \; d\mathcal H^{n-1} \\
& = \int_{B_r(0)}\Big |(F-R)R_0 \phi_{}(z)  + (Fp-b) \Big |^2 \sqrt{1+ |\nabla \phi_{}|^2(z)} \;dz
\\ & \ge  \int_{B_r(0)}\Big |(F-R)R_0  \phi_{}(z)   + (Fp-b) \Big |^2 \; dz 
\\ & = \int_{B_r(0)} \Big | \widehat{(F-R)R_0} z + g( z)  ((F-R)R_0)^{(n)} + (Fp-b)\Big |^2 \;dz 
\\ & =   \int_{B_r(0)} \Big | \widehat{(F-R)R_0} z + \tilde g( z) ((F-R)R_0)^{(n)} + c ((F-R)R_0)^{(n)}  + (Fp-b)\Big |^2 \;dz
\end{aligned} 
 $$  We let 
 $$
  A:= \widehat{(F-R)R_0},\,\,  b':= ((F-R)R_0)^{(n)}  
, \,\,  \tilde b:= c  b'  + (Fp-b).
$$ Expanding the right hand side yields 
 $$
 \begin{aligned}
e
&  =  \int_{B_r(0)} \Big | A z|^2  \; dz + |b'|^2 \int_{B_r(0)} |\tilde g (z)|^2 \;dz + |\tilde b|^2 |B_r(0)| \\ & + 2 \int_{B_r(0)}\tilde  g( z) \langle   A z ,  b' \rangle \; dz 
+  2 \int_{B_r(0)} \langle  A z ,   \tilde b \rangle  + \tilde g( z) \langle  b',  \tilde b \rangle    \; dz  
\\ & \ge  \int_{B_r(0)} \Big | A z|^2  \; dz +   |b'|^2 \int_{B_r(0)} |\tilde g (z)|^2 \;dz  + 2  \Big \langle     \int_{B_r(0)} \tilde g( z) A z \; dz  ,   b' \Big   \rangle \\ & +  2 \Big \langle  \int_{B_r(0)}  z \; dz ,   A^T \tilde b  \Big \rangle  +  2 \int_{B_r(0)} \tilde g( z) \, dz \, \langle b',  \tilde b \rangle .
\end{aligned} 
 $$ To proceed we use the facts that $\ds \int_{B_r(0)} z \, dz =0$ and   $\ds \int_{B_r(0)} \tilde g( z) \, dz  =0$ to obtain
\begin{equation}\label{e-1}
e\ge \int_{B_r(0)} \Big | A z|^2  \; dz +  |b'|^2 \int_{B_r(0)} |\tilde g (z)|^2 \;dz  + 2 \Big \langle     \int_{B_r(0)}  A z \; dz ,  \tilde g( z) b' \Big   \rangle. 
\end{equation}  

Now, we observe that there exists $0 \le \rho<1$ such that  
$$
\ds \Big |\langle     \int_{B_r(0)}  A z \; dz ,  \tilde g( z)  b' \Big   \rangle \Big | \le \rho  |b'| \Big (\int_{B_r(0)} |\tilde g (z)|^2 \;dz \Big)^\frac 12 \Big ( \int_{B_r(0)} \Big | A z|^2  \; dz \Big )^\frac 12.
$$ If either $b'=0$ or $A=0$ the claim is obvious. Otherwise if no such $\rho$ exists, by the Cauchy-Schwartz inequality  the two vector valued functions $\tilde g(z) b'$ and  $Az$ must be co-linear over $B_r(0)$ and for a $\lambda\neq 0$ we have
$$
\tilde g(z) b' = \lambda Az.
$$  Since $b'\neq 0$, this implies that $\tilde g$ is linear and hence $B(x,r)$ must be a portion of a hyperplane, in contradiction to our original choice for part (ii). So in this case as $\tilde g$ is not affine, $\rho$ depends only on the local nonlinearity of $g$, irrespective of what $A$ is, and thus can be globally measured by how far $S$ is from being a hyperplane. If otherwise $g$  is affine, under the assumptions of part (i), it is  $F$ which determines a value for $\rho$ independent of $\tilde g$. Indeed, we observe that $(F-R)R_0 e_j = Ae_j = \lambda^{-1} \tilde g(e_j) b'$ for $j=1,\cdots, n-1$, and $(F-R)R_0 e_n=b'$ by definition of $b'$, which implies that  $F-R$ is of rank 1, contrary to our assumption.

Following the above observation, \eqref{e-0} and \eqref{e-1} imply
\begin{equation*}
\begin{aligned}
~ & \int_{B_r(0)} | A z|^2  \; dz + |b'|^2 \int_{B_r(0)} |\tilde g (z)|^2 \;dz   \\ & \le \frac 1 {1-\rho}  \Big (  
\int_{B_r(0)}  | A z|^2  \; dz +  |b'|^2 \int_{B_r(0)} |\tilde g (z)|^2 \;dz  - 2 \rho  |b'| \Big (\int_{B_r(0)} |\tilde g (z)|^2 \;dz \Big)^\frac 12 \Big ( \int_{B_r(0)}  | A z|^2  \; dz \Big )^\frac 12 \Big )\\ &  \le \frac 1 {1-\rho}  e \le \frac{C}{1-\rho} E. 
\end{aligned}
\end{equation*}

Defining $G:=  A^T A$, we first observe that  
$$
 |A|^2 =   \tr(G). 
$$ On the other hand, since $G$ is symmetric nonnegative, $G= R_1 ^T D R_1$ with respectively diagonal  matrix
$D = {\rm diag}\{K_1, \cdots, K_{n-1}\}$ and orthogonal  matrix $R_1 \in SO(n-1)$, and $\tr(G) = \tr(D)$.  Hence, by symmetry,  
$$
\begin{aligned}
 \tr(G) \int_{B_r(0)} z^2_1\,  dz & = \sum_{j=1}^{n-1}    \int_{B_r(0)} K_j z^2_j \,  dz  =  \int_{B_r(0)}   \langle z , D z \rangle  \, dz   = \int_{B_r(0)}   \langle R_1 z , D R_1 z \rangle \, dz  \\ & = \int_{B_r(0)} \Big \langle z , G z \rangle \, dz =  \int_{B_r(0)} \Big | A z|^2  \; dz   \end{aligned}
$$  Therefore we obtain
$$
|A|^2 + |b'|^2 \int_{B_r(0)} |\tilde g (z)|^2 \;dz  \le \frac{C}{(1-\rho) r^{n+2}}  E.
$$  To see part (i), observe that if $S$ is no entirely flat, we can choose $p$ in such a manner that 
$$
\int_{B_r(0)} |\tilde g (z)|^2 \;dz >0. 
$$ Therefore we conclude with \eqref{fin-es} and finish the proof:
$$
|F-R|^2 = |(F-R)R_0|^2 = |A|^2 + |b'|^2 \le C(\Omega, S) E.
$$ If, $S$ is inside a hyperplane,  and $\tilde g$ is still affine, but non-zero,  the same conclusion holds, but this time with $C$ depending also on $F$ through $\rho$. Finally, if $\tilde g\equiv 0$, the best we obtain is
$$
|\widehat{(F-R)R_0}|^2 \le \frac{C(\Omega, F)}{r^{n+2}} E.
$$ But in this case $F$ is assumed not to be rank one connected to $SO(n)$, which means that for all $\tilde R\in SO(n)$, $m:= {\rm rank}((F-\tilde R)R_0)$ is either $0$ or at least $2$. If $m=0$ for some $\tilde R$, then $F=\tilde R \in SO(n)$, in which case we note that 
$$
|F-R|^2  = |(F-R)R_0|^2 = |(\tilde R- R)R_0|^2 \le n^2 |\widehat{(\tilde R - R)R_0}|^2   = n^2 |\widehat{(F-R)R_0}|^2 \le  \frac{C(\Omega, F)}{r^{n+2}} E, 
$$ establishing  \eqref{fin-es}, where we used the fact that for any rotation, 
$$
R=[v_1, \cdots, v_n] \in SO(n), 
$$ 
the last column is {\it the exterior product} 
$$
v_n= v_1 \wedge v_2 \wedge \cdots \wedge v_{n-1}
$$
of the first $n-1$ columns, which yields for any two rotations $R_1, R_2\in SO(n)$:
$$
 |R_1^{(n)} -R_2^{(n)}|  \le  (n-1) |\widehat{R_1} -\widehat{R_2}|,
$$  finally implying 
$$
|R_1- R_2|^2=  |(R_1 -R_2)^{(n)}|^2 + |\widehat{R_1 -R_2}|^2    \le n^2 |\widehat{R_1 -R_2}|^2. 
$$ If on the other hand, $m\ge 2$, it is straightforward to see that 
$$
C_F:= \inf_{\tilde R\in SO(n)} |\widehat{(F-\tilde R)R_0}|^2 >0. 
$$ We therefore obtain
$$
|F-R_F|^2 = \frac{1}{C_F}  {\rm dist}^2 (F, SO(n))  C_F  \le  \frac{C'_F}{C_F} |\widehat{(F-R)R_0}|^2 \le  \frac{C(\Omega, F)}{r^{n+2}} E,
$$ and
$$
|R_F-R|^2 \le 4n \le \frac{C}{C_F} |\widehat{(F-R)R_0}|^2 \le \frac{C(\Omega, F)}{r^{n+2}} E,
$$ which once again completes the proof. 
 
\end{proof}

\begin{cor}\label{fixed-Id-2}
Let $Q = (0,1)^3$ be  the unit cube.  For $0<r_0<1$, let $$S(x',r_0):= B(x',r_0) \cap (0,1)^2,$$ and  
$$
\mathcal{A}_{r_0}:= \{y\in W^{1,2}(Q, \R^3);\,\, \exists x' \in (0,1)^2 \,\, y|_{S(x',r_0)}(x)= x\}. 
$$   Then there exists $C=C(r_0)$ such that for all $y\in \mathcal{A}_{r_0}$, 
 $$
 \|\nabla y - {\rm Id}\|^2_{L^2 (Q)} \le C \|{\dist} (\nabla y , SO(n))\|^2_{L^2(Q)}. 
 $$ 
 \end{cor} 
 
   \subsection{Compactness for bounded sequences}
  


   \begin{lemma}\label{mollified-grady}
  Let $\Omega$, $\Omega^h$ and $K$ be as above.   Then there exist constants $h_0>0$, $C>0$, $0<c<1$, depending only on $\Omega$ and $K$, such that for all $h<h_0$ and 
$y \in \mathcal{A}^h_K$ there exists a matrix valued mapping 
  $
  \widetilde F\in W^{1,2}(\Omega, \R^{3\times 3})$,  extended trivially to $\Omega^h$, for which
 the following estimates hold true:
  \begin{equation}\label{for-F-tild}
  \ds \frac 1h \|\nabla y- \widetilde F\|^2_{L^2(\Omega^h)}  \le C E^h(y), \quad   \|\nabla \widetilde F\|^2_{L^2(\Omega)} 
    \le \frac{C}{h^2} E^h(y).
  \end{equation} Moreover, if 
  $$
  \tilde d(x'):= \left \{ \begin{array}{ll} |\widetilde F(x') - {\rm Id}| & \mbox{if} \,\, x' \in K_{ch} 
  \\ {\rm dist}(\widetilde F(x'), SO(3)) & \mbox{otherwise}
  \end{array} \right .
  $$ then
  \begin{equation}\label{for-d-tild}
 \|\tilde d\|^2_{L^2(\Omega)}  \le C E^h(y), \quad   \|\tilde d\|^2_{L^\infty(\Omega)}  \le \frac{C}{h^2} E^h(y).
  \end{equation}
   \end{lemma} 
  
  \begin{proof} The   proof closely follows \cite[Theorem 6]{FJMhier}. We will need only to make minor but careful adjustments using  Corollary \ref{fixed-Id-2}. 
    
  We cover $\Omega$ with open sets $\{U_j\}_{j=0}^N$ such that $\overline U_0 \subset \Omega$ and $\partial \Omega$, and for $j=1,\cdots, N$, $U_j \cap \Omega$ is such that for an open interval $I_j\subset \R$ and a Lipschitz function $g_j: I_j \to \R$ we have
  $$
  U_j \cap \Omega= \{x\in U_j, x_1\in I_j, x_2 > g_j(x_1)\} \quad \mbox{and}  \quad   U_j \cap \partial \Omega= \{x\in U_j, x_1\in I_j, x_2 =g_j(x_1)\}
  $$ for a suitable orthonormal coordinate system adapted to $U_j\cap \Omega$.  We also consider the flattening bi-Lipschitz change of variable 
$\Phi_j : U_j \cap \overline \Omega \to \overline{\R^2_+}$ defined by $\Phi(x_1, x_2)=   (x_1, x_2 - g_j(x_1))$.  Let also $\theta_j \in C^\infty_c(\R^2)$ be a partition of unity subject to the family $U_j$, i.e.\@ for all $j\in \{0,1,\cdots, N\}$, 
 $$
 E_j:= {\rm supp}\, \theta_j \subset U_j \quad \mbox{and} \quad \ds \sum_{j=0}^N \theta_j =1.
 $$ We let for all $x'\in \Omega$ 
 $$
 \bar F(x'):= \frac 1h \int_0^h \nabla y(x', x_3) \, dx_3.  
 $$

 {\it {\bf Step 1.} Interior local estimates:} We will first construct a matrix field $\widetilde F_0$ on $E_0 \subset U_0$ with useful  local estimates.   
 
 \medskip
We choose  $h_0$ small enough such that  ${\rm dist}(E_0, \partial U_0) >3h_0$.  For $S_{0,1}=(0,1)^2$ being the open unit square in $\R^2$, we consider a standard mollifier $\phi\in C^\infty_c(Q(0,1))$, with $\phi \ge 0$ and $\int_{S_{0,1}} \phi =1$, and we set
 $ 
 \phi_h(x):= h^{-2}\phi(x/h).
 $  and we define for each $x'\in E_0$,  
 $$
 \widetilde F_0(x'):= \phi_h \ast \bar F (x') = \frac 1{h^3} \int_{S_{x',h} \times (0,h)} \phi\Big (\frac{x'-z'}{h}\Big ) \nabla y (z)\, dz'  \, dz_3,
 $$ where $z= (z', z_3) \in \R^3$, and $S_{x',h} =   (x', x'+h)^2$ is the square of edge size $h'$ with  its lower left corner at $x'$.    
 
 \medskip

We now consider the cube $Q= \mathcal{Q}_0(x',h)= S_{x',h} \times (0,h)$ and apply Theorem \ref{rig-es-th}  to $y|_Q$.  Therefore there exists $R_{x'}\in SO(n)$ such that 
\begin{equation}\label{rig-Rx'h}
\ds \|\nabla y  -R_{x'}\|^2_{L^2( \mathcal{Q}_0(x',h))} \le C \int_{ \mathcal{Q}_0(x',h)} {\rm dist}^2 (\nabla y , SO(3)). 
\end{equation}  Moreover, fixing  a constant $0<c_0<1$, in view of Corollary \ref{fixed-Id-2}, $R_{x'}$  can be chosen to be equal to the identity matrix ${\rm Id}$ with a uniform constant $C>0$ in \eqref{rig-Rx'h} if  $K_{c_0 h} \cap \overline{S_{x',h}} \neq \emptyset$.  Indeed, in that case,  there exists $\tilde x' \in K_{c_0 h} \cap S_{x',h}$, such that setting $r_0= 1-c_0$ 
we have $y= {\rm id}$ on 
$$
B(\tilde x', r_0h) \cap S_{x',h} \subset K_{h} \cap S_{x',h} \subset \partial Q,
$$ and Corollary \ref{fixed-Id-2} applies after a proper translation and rescaling.
Also note that in this case $C$ is independent of $h$ as the boundary conditions and  the cube estimates  in both  Theorem \ref{rig-es-th} and Corollary \ref{fixed-Id-2} are invariant under dilations and translations.  Now, since 
$$
\frac 1{h^3}  \int_{ \mathcal{Q}_0(x',h)} \phi\Big (\frac{x'-z'}{h}\Big ) \, dz = 1,  
$$ letting $d\mu =  h^{-3} \phi ((x'-z')/h) \, dz$ and  applying Jensen's inequality yields
$$
|\widetilde F_0 (x') - R_{x'}|^2 = \Big | \int_{\mathcal{Q}_0(x',h)} ( \nabla y - R_{x'}) \, d\mu \Big  |^2   \le
 \int_{\mathcal{Q}_0(x',h)}  |\nabla y - R_{x'}|^2 \, d\mu,    
$$ which implies, using a uniform bound on  $\phi$ and the rigidity estimate,
\begin{equation}\label{loc-Fh}
|\widetilde F_0 (x') - R_{x'}|^2   \le \frac{C}{h^3} \int_{\mathcal{Q}_0(x',h)} {\rm dist}^2 (\nabla y , SO(3)).
\end{equation} To obtain a bound on $\nabla \widetilde F_0$ we proceed in a similar manner. For all $\tilde x'\in S_{x',h}$ we note that 
$$
\int_{\mathcal{Q}_0(x',h)} \nabla \phi\Big (\frac{\tilde x'-z'}{h}\Big ) \, dz = 0,  
$$  and we obtain  this time using Cauchy-Schwarz inequality and a uniform bound on $\nabla \phi$
 $$
 \begin{aligned}
 |\nabla \widetilde F_0(\tilde x')|^2 & =  \Big | h^{-4} \int_{\mathcal{Q}_0(x',h)} (\nabla y - R_{\tilde x',h})   \nabla \phi \Big (\frac{\tilde x'-z'}{h} \Big) \, dz  \Big |^2   \\  & \le  \Big (  \int_{\mathcal{Q}_0(x',h)} |\nabla y - R_{\tilde x' ,h}|^2 \, dz \Big )   \int_{\mathcal{Q}_0(x',h)}   \Big | h^{-4} \nabla \phi \Big (\frac{\tilde x'-z'}{h}\Big ) \Big |^2   \, dz 
  \end{aligned}
 $$  yielding 
 \begin{equation}\label{loc-gradF}
 \begin{aligned}
 \forall \tilde x' \in S_{x',h} \quad  |\nabla \widetilde F_0(\tilde x')|^2  
 & \le \frac{C}{h^5}   \int_{\mathcal{Q}_0(\tilde x',h)} {\rm dist}^2 (\nabla y , SO(3))
 \\ & \le  \frac{C}{h^5}  \int_{S_{x',2h}\times (0,h)} {\rm dist}^2 (\nabla y , SO(3)),
 \end{aligned}
 \end{equation}  and henceforth,   integrating the pointwise estimate in \eqref{loc-gradF} yields
$$
 \forall \tilde x' \in S_{x',h} \quad  |\widetilde F_0(\tilde x') - \widetilde F_0(x')|^2 \le \frac{C}{h^3}  \int_{S_{x',2h}\times (0,h)} {\rm dist}^2 (\nabla y , SO(3)), 
$$ which combined with \eqref{loc-Fh} gives
\begin{equation}\label{loc-d}
 \forall \tilde x' \in S_{x',h} \quad  |\widetilde F_0(\tilde x') - R_{x'}|^2 \le \frac{C}{h^3}  \int_{S_{x',2h}\times (0,h)} {\rm dist}^2 (\nabla y , SO(3)). 
\end{equation} On the other hand, applying the triangle inequality we have
$$
 \forall z \in \mathcal{Q}_0(x',h)  \quad |\widetilde F_0(z') -  \nabla y(z)| \le |\widetilde F_0(z') - R_{x'} |  + |R_{x'}-  \nabla y(z)|  
$$ which leads, in view of \eqref{rig-Rx'h} and \eqref{loc-d},  to 
$$
 \int_{\mathcal{Q}_0(x',h)} |\widetilde F_0(z') -  \nabla y(z)|^2 \, dz  \le  C  \int_{S_{x',2h} \times (0,h)} {\rm dist}^2 (\nabla y , SO(3)).
$$  We cover $E_0$ with a lattice of non-overlapping squares $S_{x'_i,h}$  such that $S_{x'_i,2h} \subset U_0$ for all $i$. Summing the last estimate over $i$ we obtain
\begin{equation}\label{int-F}
 \frac 1h \int_{E_0\times (0,h)} |\widetilde F_0(z') -  \nabla y(z)|^2 \, dz  \le  \frac{C}{h}  \int_{U_0 \times (0,h)} {\rm dist}^2 (\nabla y , SO(3)).
\end{equation} Similarly, intergrating  \eqref{loc-gradF} over the cubes $\mathcal{Q}_0(x'_i, h)$ and summing up over $i$ yields
\begin{equation}\label{int-gradF}
\int_{E_0} |\nabla \widetilde F_0(z')|^2 \, dz' = \frac 1h  \int_{E_0\times (0,h)} |\nabla \widetilde F_0(z')|^2 \, dz  \le  \frac{C}{h^3}  \int_{U_0 \times (0,h)} {\rm dist}^2 (\nabla y , SO(3)).
\end{equation}

For further use we establish a local version of  the $L^2$ estimate in \eqref{for-d-tild}.  For $x'\in E_0$ we define
   $$
  \tilde d_0(x'):= \left \{ \begin{array}{ll} |\widetilde F_0(x') - {\rm Id}| & \mbox{if} \,\, x' \in K_{c_0h} 
  \\ {\rm dist}(\widetilde F_0(x'), SO(3)) & \mbox{otherwise}
  \end{array} \right .
  $$  We cover  $E_0$ as before by the non-overlapping squares $S_{x'_i, h}$. If $K_{c_0h} \cap \overline{S_{x'_i, h}} \neq \emptyset$, then $R_{x'_i} = {\rm Id}$ and so integrating \eqref{loc-d} on $S_{x'_i, h}$ gives
$$
\int_{S_{x'_i, h}} d_0^2 (z') \,dx'\le  \int_{S_{x'_i, h}} |\widetilde F(z') - {\rm Id}|^2 \,dz'
\le   \frac{C}{h}  \int_{S_{x', 2h} \times (0,h)} {\rm dist}^2 (\nabla y , SO(3)).
$$ Otherwise, if $K_{c_0h} \cap \overline{S_{x'_i, h}} = \emptyset$ we have also by definition of $d_0$  and \eqref{loc-d}
$$
\begin{aligned}
\int_{S_{x'_i, h}} d_0^2 (z') \,dz' & = \int_{S_{x'_i, h}} {\rm dist}^2(\widetilde F_0(z'), SO(3))  \, dz'  
\\ & \le  \int_{S_{x'_i, h}} |\widetilde F(z') - R_{x'}|^2 \,dz'
\le   \frac{C}{h}  \int_{S_{x', 2h} \times (0,h)} {\rm dist}^2 (\nabla y , SO(3)).
\end{aligned}
$$ Summing up  the last two inequalities over $i$ we obtain
  \begin{equation}\label{d0-L2}
 \|\tilde d_0\|^2_{L^2(E_0)}  \le \frac{C}{h}  \int_{U_0 \times (0,h)} {\rm dist}^2 (\nabla y , SO(3)). 
  \end{equation}

   {\it {\bf Step 2.} Boundary  estimates:}  We will construct a matrix field $\widetilde F_j$ on $E_j \cap \Omega$ with useful  local estimates.      
  One again, letting $\xi'= \Phi(x')$, we note that  
 $$
\widetilde F_j(x') := \phi_h \ast (\bar F \circ  \Phi_j^{-1})   (\xi')
 $$  is well-defined for all $x' \in E_j \cap\Omega$ and $h$ small enough, as the square $S_{\xi', h}$ lies entirely  within the open set $\Phi_j (U_j\cap \Omega)$ in the upper half-plane. 
 \medskip
 
 In this step, we apply Theorem \ref{rig-es-th}   (respectively Corollary \ref{fixed-Id}(i)) to the Lipschitz domains 
 ${\mathcal Q}_j (x', h):= \Phi_j^{-1}(S_{\xi', h})\times (0,h)$, noting that  the constants in Theorem \ref{rig-es-th}  and  Corollary \ref{fixed-Id}(i)) are invariant under bi-Lipschitz  transformations  of  domains and  of  boundary conditions under the same transformations. Hence we have
 \begin{equation}\label{rig-Rx'h-boundary}
\ds \|\nabla y  -R_{x'}\|^2_{L^2(\mathcal{Q}_j(x',h))} \le C \int_{\mathcal{Q}_j(x',h)} {\rm dist}^2 (\nabla y , SO(3)), 
\end{equation} where
$
R_{x'} ={\rm Id} 
$ whenever $K_{c_jh} \cap \overline{\Phi_j^{-1}(S_{\xi', h})} \neq \emptyset$ for a suitable $0<c_j<1$.  Following in the footsteps of \cite[Theorem 6, Step 2]{FJMhier} the following estimates similar as in Step 1 are achieved:

\begin{equation}\label{bd-Fh}
|\widetilde F_j (x') - R_{x'}|^2   \le \frac{C}{h^3} \int_{\mathcal{Q}_j(x',h)}{\rm dist}^2 (\nabla y , SO(3)),
\end{equation}

\begin{equation}\label{bd-F}
 \frac 1h \int_{(E_j \cap \Omega)_\times (0,h)} |\widetilde F_j(z') -  \nabla y(z)|^2 \, dz  \le  \frac{C}{h}  \int_{(U_j\cap \Omega)\times (0,h)} {\rm dist}^2 (\nabla y , SO(3)),
\end{equation}  and 
\begin{equation}\label{bd-gradF}
\int_{E_j\cap \Omega} |\nabla \widetilde F_j(z')|^2 \, dz'   \le  \frac{C}{h^3}  \int_{(U_j\cap \Omega) \times (0,h)} {\rm dist}^2 (\nabla y , SO(3)).
\end{equation} Also if for $x'\in E_j\cap \Omega$
   $$
  \tilde d_j(x'):= \left \{ \begin{array}{ll} |\widetilde F_j(x') - {\rm Id}| & \mbox{if} \,\, x' \in K_{ch} 
  \\ {\rm dist}(\widetilde F_j(x'), SO(3)) & \mbox{otherwise}
  \end{array} \right .
  $$ we can prove as before for $d_0$ in \eqref{d0-L2}:
\begin{equation}\label{dj-L2}
 \|\tilde d_j\|^2_{L^2(E_j\cap \Omega)}  \le \frac{C}{h}  \int_{(U_j \cap \Omega)\times (0,h)} {\rm dist}^2 (\nabla y , SO(3)).  
  \end{equation}

 \medskip
   
   {\it {\bf Step 3.} Gluing the interior and boundary  estimates together:}
 \medskip
   
We now set  
   $
 \ds  \widetilde F := \sum_{j=1}^N \theta_j \widetilde F_j,
   $ which is well-defined on $\Omega$. We trivially extend $\theta_j$, $\widetilde F$ to $\Omega^h$. 
 We note that $\sum \nabla \theta_j =0$, and hence
 $$
\ds \nabla y -    \widetilde F   = \sum_j \theta_j  (\nabla y  - \widetilde F_j), \quad \nabla \widetilde F = \sum_j \theta_j \nabla \widetilde F_j + \sum_j \nabla \theta_j  (\widetilde F_j - \nabla y).
 $$  
Both estimates  in \eqref{for-F-tild} immediately follow from \eqref{int-F}, \eqref{int-gradF}, \eqref{bd-F} and \eqref{bd-gradF}, with constant $C$ depending on the fixed partition of unity $\{(U_j, \theta_j)\}_{j=0}^N$, i.e.\@ only on $\Omega$. 

\medskip

To establish \eqref{for-d-tild} if $c= \min\{c_j, j=0,\dots, N\}$, we obtain by the estimates \eqref{loc-Fh} and \eqref{bd-Fh}  
$$
\ds \forall x' \in K_{ch} \quad \tilde d(x')^2 =\Big |\sum_j  \theta_j (\widetilde F_j - {\rm Id}) \Big|^2 \le   C \sum_j \Big |\widetilde F_j - {\rm Id} \Big|^2  \le \frac C{h^3}    \int_{\Omega^h} {\rm dist}^2 (\nabla y , SO(3)),  
 $$  in view of the fact that if $x'\in K_{c_jh}$,  then necessarily $K_{c_jh} \cap \overline{\Phi_j^{-1}(S_{\xi',h})} \neq \emptyset$ ($K_{c_0h} \cap \overline{S_{x',h}} \neq \emptyset$ for $j=0$) and thus
  the rotation 
  $
 R_{x'} = {\rm Id},  
  $   as previously established. 
 
 Otherwise, for any $x'\in \Omega$,  let 
 $$
 J(x'):=  \{j\in \{0, \dots, N\}; \,\, \theta_j (x') \neq 0\}. 
 $$ Note that $j \in J(x')$ implies $x'\in E_j$. Since the $\Phi_j$ are bi-Lipschitz, there exists a constant $C_0>0$ depending on the domain $\Omega$ only such that 
 $$
 \bigcup_{j\in J(x')} \mathcal Q_j (x', h) \subset \ddot{\mathcal{Q}}:=  \mathcal (B(x', C_0h) \cap \Omega) \times (0,h),  
 $$ where $B(x', r)$ is the disk of radius $r$ centered at $x'$. Applying once again Theorem \ref{rig-es-th} to $y$ on $\ddot{\mathcal{Q}}$, and noting that the constant $C$ still depends only on the Lipschitz constant of $\partial \Omega$, we have for some rotation $\ddot R_{x'}\in SO(3)$
\begin{equation}\label{rig-R-biggercube}
\ds \|\nabla y  -\ddot R_{x'}\|^2_{L^2(\ddot{\mathcal{Q}})} \le C \int_{ \ddot{\mathcal{Q}}} {\rm dist}^2 (\nabla y , SO(3)). 
\end{equation}

We note that, since for all $j\in J(x')$, $\mathcal{Q}_j(x',h) \subset \ddot{\mathcal{Q}}$, we can use $\ddot{R}_{x'}$ instead of the rotations used in  \eqref{loc-Fh} and \eqref{bd-Fh}, this time obtaining the new bounds
$$
 \forall j \in J(x')   \,\, \quad  |\widetilde F_j (x')- \ddot R_{x'}  | \le  \frac{C}{h^3} \int_{\mathcal{Q}_j(x',h)}| \nabla y - R_{x'}|^2  \le  \frac {C}{h^3} \int_{\ddot{\mathcal{Q}}} {\rm dist}^2 (\nabla y , SO(3)) \le  \frac {C}{h^2} E^h (y)
$$ which establishes
$$
{\rm dist}^2 (\widetilde F(x'), SO(3)) \le |\widetilde F(x') - R_{x'}|^2 =  \Big |\sum_{j\in J(x')} \theta_j (x') (\widetilde F_j (x') - \ddot R_{x'})\Big |^2  \le  \frac {C}{h^2} E^h(y)
$$ as required for completing the $L^\infty$ estimate for $\tilde d$ in  \eqref{for-d-tild} when $x'\notin K_{ch}$.

\medskip

To complete the proof of  \eqref{for-d-tild}, it remain to prove the $L^2$ estimate for $\tilde d$. We first define for all $x= (x', x_3) \in \Omega^h$ 
   $$
  d(x):= \left \{ \begin{array}{ll} |\nabla y (x) - {\rm Id}| & \mbox{if} \,\, x' \in K_{ch} 
  \\ {\rm dist}(\nabla y (x), SO(3)) & \mbox{otherwise}
  \end{array} \right .
  $$ and we note that  if $x = (x', x_3) \in (E_j \cap \Omega ) \times (0,h)$, we have  
  $$
  d(x) \le|  \nabla y(x) -\widetilde F_j(x') | +   \tilde d_j(x').
  $$ Hence the above $L^2$ estimates \eqref{d0-L2} and \eqref{dj-L2} on $\tilde d_j$ obtained  in Steps 1 and 2, alongside \eqref{int-F} and \eqref{bd-F} imply
  $$
\frac 1h \int_{(E_j \cap \Omega) \times (0,h)} \nabla^2 \tilde d_j(x)\, dx 
\le \frac Ch \int_{(U_j \cap \Omega) \times (0,h)} {\rm dist}^2 (\nabla y(x), SO(3)) \,dx.
  $$ Summing over $j$ gives
  $$
  \ds \frac 1h \|d\|^2_{L^2(\Omega^h)} \le C E^h(y), 
  $$ which combined with \eqref{for-F-tild} proves the $L^2$ estimate in \eqref{for-d-tild} in view of the fact that 
$$
\forall x \in \Omega\quad  \ds \tilde d(x') \le |F(x') - \nabla y(x)|  + d(x). 
$$  \end{proof}

 \begin{lemma}\label{tildF'}
 Let $\Omega$, $\Omega^h$, $K$ be as defined above. Assume that $h_0>0$, $0<c<1$ be as in Lemma \ref{mollified-grady} and set $\bar c=c/2$. Then, given $h<h_0$, $y\in \mathcal{A}^h_K$, there exists a  matrix field $\widetilde F'$  which is equal to the identity matrix {\rm Id} on $K_{\bar ch}$ such that the estimates \eqref{for-F-tild} and \eqref{for-d-tild} still hold true (possibly with a new constant $C$) for $\widetilde F'$, i.e.\@
   \begin{equation}\label{for-F-tild'}
  \ds \frac 1h \|\nabla y- \widetilde F'\|^2_{L^2(\Omega^h)}  \le C E^h(y), \quad   \|\nabla \widetilde F'\|^2_{L^2(\Omega)} 
    \le \frac{C}{h^2} E^h(y), 
  \end{equation} and 
  \begin{equation}\label{F'-infty}
  \|{\rm dist}^2 (\widetilde F, SO(3))\|^2_{L^\infty(\Omega)} \le \frac C{h^2} E^h(y). 
\end{equation}
 \end{lemma}
 
 \begin{proof}
 
 We let $\widetilde F: \Omega \to \R^{3\times 3}$ be as in the statement of Lemma  \ref{mollified-grady}. 
 We introduce a cut-off function $\psi: \R^2 \to \R$
 $$
\forall x \in \R^2  \quad 
 \psi (x)= \left 
 \{ \begin{array}{lll} 
 \ds \frac{{\rm dist} (x, K_{\bar ch}) }{\bar ch}   &  \mbox{if} \,\, {\rm dist} (x, K) \le ch 
 \\  1 &   \mbox{otherwise}.
  \end{array} \right .
 $$ It can be shown that $0 \le \psi \le 1$, ${\rm supp}\, \psi = \overline {K_{\bar ch}}$, ${\rm supp}\, (1-\psi) = \R^2 \setminus  {K_{ch}}$, and that $\psi$ is a Lipschitz function with 
\begin{equation}\label{psi-lip}
\|\nabla \psi\|_{\infty} \le \frac 1{\bar ch}. 
 \end{equation}
  We now define
  $$
  \widetilde F' = \psi \widetilde F + (1-\psi) \rm {Id}. 
  $$ Note that  $\widetilde F' = {\rm Id}$ on $K_{\bar ch}$ as required. Also $\widetilde F' = \widetilde F$ on $\Omega \setminus K_{ch}$. 
   We have
   $$
   |\widetilde F' - \widetilde F| = | \psi \widetilde F  + (1-\psi)  {\rm Id} - (\psi \widetilde F + (1-\psi)  \widetilde F |  = 
   |(1-\psi) (\widetilde F - {\rm Id})| \le \tilde d.   
    $$   Hence  
    $$
    {\rm dist}(\widetilde F' , SO(3)) \le  |\widetilde F' - \widetilde F|  +  {\rm dist}(\widetilde F, SO(3))  \le 2\tilde d \le \frac{C}{h^2} E^h(y) 
    $$ as required for \eqref{F'-infty}. 
 Combining the first estimates in  \eqref{for-F-tild} and \eqref{for-d-tild} with  $|\widetilde F' - \widetilde F| \le \tilde d$ we  also obtain  
    $$
    \frac 1h \|\nabla y - \widetilde F'\|_{L^2(\Omega^h)} \le C E^h(y).
    $$
  Now we have
    $$
    \nabla \widetilde F' = \nabla \psi \otimes (\widetilde F  - {\rm Id})+ \psi \nabla \widetilde F , 
    $$ which yields, in view of  \eqref{psi-lip} and the fact that $\nabla \psi =0$ on the complement of $K_{ch}$,
    $$
    | \nabla \widetilde F'| \le \frac Ch \tilde d + |\nabla \widetilde F|.
    $$ Once again, the $L^2$ bounds on $\nabla \widetilde F$ and $\widetilde d$ in \eqref{for-F-tild} and \eqref{for-d-tild} imply
     $$
    \frac 1h \|\nabla \widetilde F'\|_{L^2(\Omega^h)} \le \frac{C}{h^2} E^h(y),
    $$ completing the proof of  \eqref{for-F-tild'}.
     \end{proof}
   \begin{rem}
In general, by applying a projection on a large ball containing $SO(3)$, we can also assume that $|\widetilde F'| \le C$. By \eqref{F'-infty}, this will not be needed under the assumption $E^h \le Ch^2$.  
  \end{rem}
  
In the following statement we now arrange $\tilde F'$ in such a way that it takes values only into $SO(3)$.
 
 \begin{cor}    \label{for-R-cor}
 Let $\Omega$, $\Omega^h$, $K$ be as defined above. Then there exist constants $h_0>0$, $C>0$, $0<\bar c<1$, $\delta_0>0$, depending only on $\Omega$ and $K$, such that if for any $h<h_0$, $y \in \mathcal{A}^h_K$ we have $E^h(y)\le \delta_0 h^2$, then there exists a matrix valued mapping 
 $
 R\in W^{1,2}(\Omega, SO(3))$ such that  $R|_{K_{\bar ch}} \equiv {\rm Id}
$ and  the following estimates hold true:
  \begin{equation}\label{for-R}
  \ds \frac 1h \|\nabla y- R\|^2_{L^2(\Omega^h)}  \le C E^h(y), \quad   \|\nabla R\|^2_{L^2(\Omega)} 
    \le \frac{C}{h^2} E^h(y).
\end{equation}
 \end{cor}
 
 \begin{proof}
We let $\widetilde F$, $h_0, \bar c$ be chosen according to Lemma \ref{for-F-tild'}. The result follows similarly as in \cite[Remark 5]{FJMhier}: Under the assumptions of the Corollary \ref{for-R-cor}, \eqref{F'-infty} implies that $\widetilde F'$ is in a $(C\delta_0)^{\frac12}$-tubular neighborhood $\mathcal U$ of $SO(3)$.  For $\delta_0$ small enough, the orthogonal projection $\pi$ from $\mathcal U$  onto $SO(3)$ is well-defined and Lipschitz. We let $R:= \pi \circ \widetilde F'$, and note that
$$
\begin{aligned}
|\nabla y - R| & \le |\nabla y - \widetilde F'|  + |\widetilde F'- R|\\ &  = |\nabla y - \widetilde F'|  + {\rm dist}(\widetilde F', SO(3))
\\ & \le  2|\nabla y - \widetilde F'|   + {\rm dist}(\nabla y, SO(3)),
\end{aligned}
$$ and that 
$$
|\nabla R| \le \|\nabla \pi\|_{L^\infty} |\nabla \widetilde F'|,
$$
which combined with \eqref{for-F-tild'} imply  together \eqref{for-R}. 
 \end{proof}

 Let for $x_3 \in (0,1)$, $\tilde y^h (x', x_3):= y^h (x', hx_3) = ((y^h)'(x', hx_3), y^h_3 (x', hx_3))$. We consider the mappings 
 $$
  {\rm id}^h : \Omega^h \to \R^3, \quad {\rm id}^h (x',x_3) = (x',x_3); \quad {\rm id} : \Omega^1 \to \R^3 \quad {\rm id} (x',x_3) = ({\rm id}^h)' (x', x_3) = (x',0),
 $$  and the displacement fields
 \begin{equation}\label{def-uh}
 u_h\in  W^{1,2} (\Omega,\R), \quad  U_h:= \frac 1h \int_0^h (y^h_3 (\cdot, x_3) - hx_3)   \, dx_3,\quad  u_h:= \frac{1}{h^{\beta/2-1}} U_h,
 \end{equation}
\begin{equation}\label{def-wh}
w_h \in W^{1,2} (\Omega, \R^2), \quad W_h:= \frac 1h \int_0^h ((y^h)' - {\rm id})  (\cdot, x_3)  \, dx_3 , \quad w_h:= \frac{1}{h^{\delta}}W_h \,\, \mbox{for} \,\, \delta:= \min\{\beta-2, \beta/2\}.
\end{equation}
\begin{theorem}\label{compactness}
 Let $\Omega$, $\Omega^h$, $K$ be as above, and assume that ${\rm Cap}_{1,2}(K)>0$ and $\beta>2$. Assume that for a sequence as $h\to 0$, $y^h \in \mathcal{A}^h_K$, $E^h(y^h) \precsim h^\beta$.    Then, up to a subsequence as $h\to 0$, we have
 \begin{itemize}
 \item[(i)] $\ds \frac 1h \|y^h  - {\rm id}^h\|^2_{W^{1,2}(\Omega^h)}\longrightarrow 0$, $($in particular $\tilde y^h  \longrightarrow  {\rm id}$  in $W^{1,2}(\Omega^1))$.
 \item[(ii)] $\ds u_h  \xrightarrow{\mbox{\tiny strongly  in} \,\, W^{1,2}} u \in H^{2}_{0,K}(\Omega)$.
  \item[(iii)]  $\ds w_h \xrightarrow{\mbox{\tiny weakly in} \,\, W^{1,2}} w \in W^{1,2}(\Omega)$.  
 \item[(iv)] If $2<\beta<4$,  then $\ds w_h \xrightarrow {\mbox{\tiny strongly in} \,\, W^{1,2}} w \in H^1_{0,K}(\Omega)$; and 
 $$
 \ds {\rm sym}\, \nabla w + \frac 12 \nabla u \otimes \nabla u \equiv 0\quad  \mbox{in} \,\, \Omega. 
  $$ 
   
{\item[(v)] For $\beta\ge 4$,  assume that $\ds {\rm sym} \nabla w_h$ converges strongly in $L^2(\Omega)$, and that moreover the family
 $$
\mathcal{Y}:= \frac 1 {h^\beta}  \Big ( \frac 1h  \int_0^h {\dist}^2(\nabla y^h (\cdot, x_3), SO(3))\;dx_3  \Big)
 $$  is equi-integrable over $\Omega$, in the sense that for all $\varepsilon>0$, there exists $\delta>0$ such that for all $B\subset \Omega$, $|B|<\delta$, we have  
 $$
\frac 1 {h^\beta} \Big ( \frac 1h  \int_{B \times (0,h)} {\dist}^2(\nabla y^h, SO(3))\;dx \Big) < \varepsilon.
 $$ Then $w_h \xrightarrow{\mbox{\tiny strongly in} \,\, W^{1,2}} w\in H^1_{0,K}(\Omega)$. }
 
 \end{itemize}
 \end{theorem}
   
\medskip

  \begin{proof}

  Let $h_0, \delta_0, \bar c, C$ be as in Corollary \ref{for-R-cor}. By the assumption $\beta>2$, and if necessary by choosing a smaller value for $h_0$, we can assure that $h<h_0$ implies $E^h(y^h) \le \delta_0 h^2$ and so  Corollary \ref{for-R-cor} applies: For each $y^h$ with $h<h_0$, we denote the associated rotation field by $R^h$, and hence \eqref{for-R} implies
    \begin{equation}\label{for-Rh}
  \ds \frac 1h \|\nabla y^h- R^h\|^2_{L^2(\Omega^h)}  \le C  h^\beta, \quad   \|\nabla (R^h- {\rm Id})\|^2_{L^2(\Omega)} 
    \le Ch^{\beta-2},
\end{equation}
where $R^h \equiv {\rm Id}$ on $K_{\bar c h}$. We  let 
\begin{equation}\label{for-Ah}
A^h:= \frac{1}{h^{\beta/2-1}} (R^h -{\rm Id}). 
\end{equation} By the second estimate in   \eqref{for-Rh}, $\nabla A^h$ is uniformly bounded in $L^2(\Omega)$ and that $A^h \equiv 0$ on $K_{\bar ch}$. Moreover, by the Poincar\'e inequality proved in Corollary \ref{qe-poinc} the sequence $A^h$ satisfies 
$$
\|A^h  \|_{L^2(\Omega)} \le C \|\nabla A^h\|_{L^2(\Omega)} 
$$ and hence up to a subsequence, $A^h$ converges weakly in $W^{1,2}$ to $A$.    Note that by Lemma \ref{QeLemma0}, we have 
\begin{equation}\label{AonK}
A\in H^1_{0,K}(\Omega).
\end{equation} We have moreover
\begin{equation}
\label{symRh}
{\rm sym}\, R^h -{\rm Id} = -\frac 12  (R^h - {\rm Id})^T (R^h - {\rm Id}), 
\end{equation} which implies that 
$$
{\rm sym} \, A^h  =  -\frac 12    h^{\beta/2-1}  (A^h)^T A^h, 
$$ and as a consequence, passing to the limit, we obtain that ${\rm sym}\, A \equiv 0$ in $\Omega$.  Now,  using
\eqref{symRh}, this time rescaled with $h^{\beta-2}$, Sobolev embeddings, and the fact that $A^T A = -A^2$ for any skew-symmetric matrix, we obtain  that 
\begin{equation}\label{symRh-2}
\forall p<2\quad B^h:= \frac{1}{h^{\beta-2}} ({\rm sym}\, R^h -{\rm Id})  = - \frac 12 (A^h)^T A^h \rightharpoonup \frac{A^2}{2} \,\,\mbox{weakly in} \,\, W^{1,p}(\Omega),
\end{equation}  which also implies the strong convergence of $B^h$  in $L^q(\Omega)$ for all $q<
\infty$. 
\medskip
 
 We first note that the convergence of $u_h$  and $w_h$ in $W^{1,2}$ as we shall prove below in parts (ii) and (iii) imply (i) in a straightforward manner. Indeed,  we first remark that  by compactness of Sobolev embeddings, $A^h$ must converge strongly to $A$ in $L^2(\Omega)$.  Now  the first estimate in  \eqref{for-Rh} implies
\begin{equation}\label{grad-to-A}
\begin{aligned}
\ds \frac 1h   \Big \|\frac{1}{h^{\frac{\beta}{2}-1}} (\nabla y^h- {\rm Id}) - A \Big \|^2_{L^2(\Omega^h)} 
&\le  \frac Ch  \Big (  \Big \|\ \frac{1}{h^{\frac{\beta}{2}-1}} (\nabla y^h- R^h)\Big \|^2_{L^2(\Omega^h)}  
 +   \|A^h- A \|^2_{L^2(\Omega^h)} \Big )  
\\ & \le  C(h^2 + \|A^h - A\|^2_{L^2(\Omega)})  \xrightarrow{h \to 0} 0,
\end{aligned}
\end{equation}  which in particular yields
\begin{equation}\label{grad-to-0}
\frac 1h \|\nabla (y^h - {\rm id}^h)\|^2_{L^2(\Omega^h)} \le Ch^{\beta-2}  \xrightarrow{h \to 0} 0.
\end{equation}
To obtain an $L^2$ estimate on $y^h - {\rm id}^h$, for a.e.~$x'\in \Omega$ we apply the Poincar\'e inequality on the interval $\{x'\} \times (0,h)$ and integrate over $\Omega$ to obtain  
$$
\begin{aligned}
\frac 1h \Big \|(y^h)'-({\rm id}^h)' - W_h\Big \|^2_{L^2(\Omega^h)} 
&  = \frac 1h  \Big \|(y^h)'-{\rm id} -  \fint_0^h ((y^h)'-{\rm id}) (\cdot, x_3) \, dx_3    \Big \|^2_{L^2(\Omega^h)} 
\\ & \le Ch^2 \Big (\frac 1h \|\partial_3 ((y^h)'-({\rm id}^h)') \|^2_{L^2(\Omega)} \Big )
\\ & \le Ch^2 \Big (\frac 1h \|\nabla (y^h -{\rm id}^h) \|^2_{L^2(\Omega)} \Big )
\end{aligned}
$$ for the first two components of $y^h - {\rm id}^h$, and 
$$
\begin{aligned}
\frac 1h \Big \|y_3^h- {\rm id}^h_3 - U_h\Big \|^2_{L^2(\Omega^h)} 
&  = \frac 1h  \Big \|y_3^h - {\rm id}^h_3  -  \fint_0^h (y_3^h  - {\rm id}^h_3) (\cdot, x_3)  \, dx_3    \Big \|^2_{L^2(\Omega^h)} 
\\ & \le Ch^2 \Big (\frac 1h \|\partial_3 (y_3^h - {\rm id}^h_3) \|^2_{L^2(\Omega)} \Big )
\\ & \le Ch^2 \Big (\frac 1h \|\nabla (y^h -{\rm id}^h) \|^2_{L^2(\Omega)} \Big ).
\end{aligned}
$$ Combining the last two estimates with \eqref{grad-to-0} we hence obtain
$$
\begin{aligned}
\frac 1h  \|y^h - {\rm id}^h\|^2_{L^2(\Omega^h)} 
& \le C (\|U_h\|^2_{L^2(\Omega)} + \|W_h\|^2_{L^2(\Omega)} + h^{\beta}) 
\\ & \le C( h^{\beta-2} \|u_h\|^2_{L^2(\Omega)}  +h^{2\delta} \|w_h\|^2_{L^2(\Omega)} + h^{\beta}) \xrightarrow{h\to 0} 0,
 \end{aligned}
 $$ where we used (ii) and (iii). 
 
\medskip 
  To see (ii), first observe that for $\nabla '= \nabla_{x'}$
 $$
 \begin{aligned}
 \|\nabla u_h - (A_{31}, A_{32})\|_{L^2(\Omega)}  
 & =  \Big \|\frac{1}{h^{\beta/2-1}} \nabla' \int_0^1   \tilde y^h_3  - hx_3 \,dx_3  - (A_{31}, A_{32}) \Big \|_{L^2(\Omega)} \\ &   =  \Big \| \int_0^1 \frac{1}{h^{\beta/2-1}} \nabla'    \tilde y_3^h  - hx_3 \,dx_3  - (A_{31}, A_{32}) \Big \|_{L^2(\Omega)} 
\\ & \le\frac {1}{\sqrt{h}} \Big \|  \frac{1}{h^{\beta/2-1}} \nabla  (y^h -{\rm id}^h) - A  \Big \|_{L^2(\Omega^h)} 
 \\ & = \frac  {1}{\sqrt{h}}  \Big \|\frac{1}{h^{\beta/2-1}}( \nabla y^h  - {\rm Id})  - A \Big \|_{L^2(\Omega^h)}  
\\ &  \le \frac  {1}{\sqrt{h}} \Big \|\frac{1}{h^{\beta/2-1}}( \nabla y^h  -  R^h)\Big \|_{L^2(\Omega^h)}  + \Big \|\frac{1}{h^{\beta/2-1}} (R^h -{\rm Id}) - A \Big \|_{L^2(\Omega)} 
 \\ &  \le C \sqrt{h}+ \|A^h   - A  \|_{L^2(\Omega)} \xrightarrow {h\to 0} 0,
   \end{aligned} 
  $$ up to a subsequence, where we used \eqref{for-Rh} and the convergence of $A^h$ to $A$ in the last line.  Therefore for $a_h:= \fint_\Omega u_h$ we obtain by the Poincar\'e-Sobolev inequalities that for some $\tilde u \in W^{2,2}(\Omega)$, with 
  $\nabla \tilde u = (A_{31}, A_{32})$,  
  $$
 \|u_h- a_h - \tilde u\|_{W^{1,2}(\Omega)} \xrightarrow{h\to 0} 0.
$$ (To see this, we first establish the weak convergence $u_h - a_h$ and then use the strong convergence of $\nabla u_h$.) Note that by Sobolev embedding theorems $\tilde u\in C^0(\overline \Omega)$, and hence the value of $\tilde u$ for all points on $K$ is well-defined. By Lemma \ref{avg-to-0} proven below, we have for q.e.~$x\in K$
$$
\begin{aligned}
\lim_{h\to 0} \fint_{B(x,\bar c h)\cap \Omega} |u_h - a_h - \tilde u|  =0 \implies \lim_{h\to 0} \fint_{B(x,\bar c h)\cap \Omega} (u_h - a_h - \tilde u)  =0, 
  \end{aligned} 
$$ which implies 
 \begin{equation}\label{uh-ah}
 \lim_{h \to 0} \fint_{B(x,\bar c h)\cap \Omega} u_h - a_h  =  \lim_{h\to 0} \fint_{B(x,\bar c h)\cap \Omega} \tilde u =\tilde u(x). 
 \end{equation} On the other hand, since $y^h= {\rm id}$ on $K_h$: 
 $$
 \begin{aligned}
 \int_{B(x,\bar c h) \cap \Omega} |U_h|^2  & \le\frac 1h   \int_{(B(x,\bar c h) \cap \Omega) \times (0,h)} |y^h_3 -  x_3|^2 
 \\ & \le Ch^2 \Big (  \frac 1h \int_{(B(x, \bar c h) \cap \Omega) \times (0,h)} |\partial_3 (y^h_3 -  x_3)|^2 \Big ) 
  \\ & \le Ch^2 \Big (  \frac 1h \int_{(B(x, \bar c h) \cap \Omega) \times (0,h)} |\nabla y^h-   {\rm Id}|^2 \Big )
 \\ & \le Ch^2 \Big (  \frac 1h \int_{(B(x, \bar c h) \cap \Omega) \times (0,h)} |\nabla y^h -   R^h|^2 \Big ) 
 \le Ch^{\beta+2}, 
  \end{aligned}
  $$  since $R^h= {\rm Id}$ on $K_{\bar c h}$. Therefore,  since $|B(x,\bar c h) \cap \Omega| \ge C(\Omega) h^2$ we obtain
  $$
   \fint_{B(x,\bar c h) \cap \Omega} |u_h|^2 \le Ch^2, 
   $$ which implies 
   $$
   \lim_{h\to 0}    \fint_{B(x,\bar c h) \cap \Omega} |u_h| =0, 
   $$ yielding, combined with \eqref{uh-ah},  that $a=\ds  \lim_{h\to 0} a_h \in \R$ exists, and that for q.e.~$x\in K$, $\tilde u(x) =-a$.  Therefore for $u:= \tilde u+a$: 
   \[
 u_h = (u_h - a_h) + a_h \xrightarrow{\mbox{in} \,\, W^{1,2}} \tilde u + a = u \in   W^{2,2}(\Omega), 
  \] and the continuous representative of $u$ vanishes on $K$.  More precisely, we obtain by Lemma \ref{avg-to-0}, since $u_h -u$ converges strongly to $0$ in $W^{1,2}(\Omega)$:
 $$
\fint_{B(x,\bar c h) \cap \Omega} |u|  \le  \fint_{B(x,\bar c h) \cap \Omega} |u- u_h|  +   \fint_{B(x,\bar c h) \cap \Omega} |u_h|  \xrightarrow{h\to0} 0. 
$$ 
Therefore, by definition, $u\in H^1_{0,K}(\Omega)$, which alongside the fact that by \eqref{AonK} $\nabla u = (A_{31}, A_{32}) \in H^1_{0,K}(\Omega)$, completes the proof of (ii).

   \medskip
   
  To establish (iii) and (iv) we proceed in the same manner as in (ii), with a caveat.  Let $\aoverbrace[L1R]{B}$ be the $2\times 2$ upper-left sub-matrix of $B:=A^2/2$. We have
  $$
 \begin{aligned}
 \|{\rm sym}\, \nabla w_h   -  h^{(\beta-2)-\delta} & \aoverbrace[L1R] B \|_{L^2(\Omega)}  
  =  \Big \|{\rm sym}\,  \Big ( \frac{1}{h^\delta} \nabla' \int_0^1   (\tilde y^h)'   - x' \,dx_3  \Big ) -  h^{(\beta-2)-\delta} \aoverbrace[L1R] B \Big \|_{L^2(\Omega)} \\ &   =  \Big \| \int_0^1 \frac{1}{h^{\delta}} {\rm sym}  \Big ( \nabla'    ((\tilde y^h)'  - x') \,dx_3 \Big )  
  - h^{(\beta-2)-\delta}   \aoverbrace[L1R] B \Big \|_{L^2(\Omega)} 
\\ & \le\frac 1{\sqrt{h}}  \Big \|  \frac{1}{h^{\delta}} {\rm sym}  (  \nabla  (y^h -{\rm id}^h) )  - h^{(\beta-2)-\delta}  B   \Big \|_{L^2(\Omega^h)} 
 \\ & = \frac  1{\sqrt{h}}  \Big \|\frac{1}{h^\delta} {\rm sym}  ( \nabla y^h  - {\rm Id}   ) - h^{(\beta-2)-\delta}  B \Big \|_{L^2(\Omega^h)}  
\\ &  \le \frac  1{\sqrt{h}}  \Big \|\frac{1}{h^{\delta}}{\rm sym} (   \nabla y^h  -  R^h )\Big \|_{L^2(\Omega^h)}  +    \Big \|\frac{1}{h^\delta} {\rm sym}\,  (R^h -{\rm Id}) -  h^{(\beta-2)-\delta} B \Big \|_{L^2(\Omega)} 
 \\ &  \le Ch^{(\beta/2)- \delta} + h^{(\beta-2)-\delta}\|B^h   - B  \|_{L^2(\Omega)}.
   \end{aligned} 
  $$  by \eqref{for-Rh}, \eqref{symRh-2}.  If $2<\beta <  4$, then $\delta = \beta-2 <\beta/2$ and we deduce that  ${\rm sym}\, \nabla w_h$ converges in $L^2(\Omega)$ to $\aoverbrace[L1R]{B}$ by \eqref{symRh-2}.  Otherwise, if $\beta\ge 4$, then $\delta= \beta/2$ and $h^{(\beta-2)-\delta} \to 0$, and as a consequence $\|{\rm sym}\, \nabla w_h  \|_{L^2(\Omega)}$  is merely bounded.  This implies, through applying Korn's inequality on $\Omega$, that  if for all $x\in \Omega$ we set   
  $
      \tilde w_h(x) := w_h(x) - D^h x - b_h,
      $  where
      $$
      b_h:= \fint_\Omega w_h \quad \mbox{and} 
      \quad D^h = \left [ \begin{array}{cc} 0 & d_h \\ -d_h & 0 \end{array} \right ] := \fint_\Omega  {\rm skew} (\nabla w_h) 
      $$  there exists $\tilde w \in W^{1,2}(\Omega, \R^2)$, with ${\rm sym} \, \nabla \tilde w = \aoverbrace[L1R]{B} $ if $2<\beta<4$,  so that up to a subsequence
       $$
\begin{array}{ll}
\mbox{if} \,\, 2 <\beta < 4 &  \tilde w_h  \xrightarrow{in \,\, W^{1,2} } \tilde w,  \\ 
 \mbox{if}\,\, \beta\ge  4  &  \tilde w_h   \xrightarrow{weakly \,\, in \,\, W^{1,2}} \tilde w.
\end{array}
$$  If $2<\beta< 4$ we obtain through Lemma \ref{avg-to-0} that for q.e.~$x\in \overline \Omega$
$$
\lim_{h\to 0} \fint_{B(x, \bar c h)\cap \Omega} \tilde w_h  = \lim_{h\to 0} \fint_{B(x, \bar c h) \cap \Omega}  \tilde w = \tilde w^\ast (x) 
$$ by applying an extension operator of Appendix \ref{qeAppx} and \cite[Theorem 4.8.1]{EG}, where $\tilde w^\ast (x)$ is the precise representative of $\mathcal{E}(\tilde w)$. We now observe  that since $y^h= {\rm id}$ on $K_h$ and $R^h \equiv {\rm Id}$ on $K_{\bar ch}$: 
 $$ 
 \begin{aligned}
 \int_{B(x,\bar c h) \cap \Omega} |W_h|^2  & \le\frac 1h   \int_{(B(x,\bar c h) \cap \Omega) \times (0,h)} |(y^h)' (x)-  x'|^2 
 \\ & \le Ch^2 \Big (  \frac 1h \int_{(B(x, \bar c h) \cap \Omega) \times (0,h)} |\partial_3 ((y^h)' -  x')|^2 \Big ) 
  \\ & \le Ch^2 \Big (  \frac 1h \int_{(B(x, \bar c h) \cap \Omega) \times (0,h)} |\nabla y^h -   {\rm Id}|^2 \Big )
 \\ & \le Ch^2 \Big (  \frac 1h \int_{(B(x, \bar c h) \cap \Omega) \times (0,h)} |\nabla y^h -   R^h|^2 \Big ) 
 \le Ch^{\beta+2}.
  \end{aligned}
  $$    Thus,  once again through $|B(x,\bar c h) \cap \Omega| \ge C(\Omega) h^2$ we obtain
  $$
   \fint_{B(x,\bar c h) \cap \Omega} |w_h|^2 \le Ch^{\beta-2\delta}, 
   $$ which implies  for $2<\beta< 4 $: 
   $$
   \lim_{h\to 0}    \fint_{B(x,\bar c h) \cap \Omega} |w_h| =0.
   $$ Therefore we obtain for q.e.~$x\in K$
$$
 \lim_{h\to 0} \Big ( d_h  \fint_{B(x, \bar c h)\cap \Omega} y^{\perp}  \, dy +     b_h \Big )  = \lim_{h\to 0} \Big ( D^h  \fint_{B(x, \bar c h)\cap \Omega} y' \, dy +     b_h \Big ) = -\tilde w^\ast  (x).  
$$    But by regularity of the mapping $y \to y^\perp$, for all $x\in \overline \Omega$
$$
 \lim_{h\to 0}  \fint_{B(x, \bar c h)\cap \Omega} y^{\perp}  \, dy = x^\perp, 
$$  so we conclude with
$$
\mbox{for q.e.}\,\, x\in K \quad \lim_{h\to 0} (d_h x^\perp + b_h) = - \tilde w^\ast  (x).
$$ Note that the q.e.~existence of the above limit implies its existence  for at least two distinct values of $x\in K$, from which immediately follows that 
$d:= \lim_{h \to 0} d_h$  and $b:= \lim_{h\to 0} b_h$ exist and subsequently that for q.e.~$x\in K$, 
$$
 \tilde w^\ast (x) =  -(dx^\perp  + b).
$$  Now, $w_h - \tilde w_h (x) = D^h x + b_h \to d x^\perp + b $ which yields that  
$$
w_h = (w_h - \tilde w_h) + \tilde w_h \xrightarrow{in\,\, W^{1,2}} w:= \tilde w + d x^\perp + b. 
$$ Applying  Lemma \ref{avg-to-0} once again  in an argument similar to what we presented for $u$ in part (ii), yields $w \in H^1_{0,K}(\Omega)$.
 
\medskip

Finally, remember that $\aoverbrace[L1R]{F}$ is the $2\times 2$ upper-left sub-matrix of $F\in R^{3\times 3}$. First we observe that  by \eqref{grad-to-A}
$$
\begin{aligned}
\|\frac{h^\delta}{h^{\beta/2-1}}\nabla w_h - \aoverbrace[L1R] {A} \|^2_{L^2(\Omega)} 
& =  \Big \|\frac{1}{h^{\beta/2-1}}\Big (     \fint_{0}^h \nabla' ((y^h)' - ({\rm id}^h)')(\cdot, x_3) \, dx_3\Big) - \aoverbrace[L1R]{A}  \Big \|^2_{L^2(\Omega)}  
\\& =  \Big \|  \fint_{0}^h \Big (\aoverbrace[L1R]{ \frac{1}{h^{\beta/2-1}}  (\nabla y^h - {\rm Id}) - A} \Big) \,dx_3   \Big \|^2_{L^2(\Omega)}  
\\ & \le   \frac 1h \Big ( \Big \|\frac{1}{h^{\beta/2-1}}(\nabla y^h - {\rm Id}) -  A \Big \|^2_{L^2(\Omega^h)}    \Big ) \xrightarrow{h\to 0} 0.
 \end{aligned}
 $$ But since $\beta>2$, we have $\delta > \beta/2-1$, and hence, in view of the boundedness of $\|\nabla w_h\|_{L^2}$ and by passing to the limit in $h\to 0$ we obtain $\aoverbrace[LR]{A}=0$. Remember that $A^T = -A$, and that $\nabla  u  = (A_{31}, A_{32})$. Hence, overall, straightforward calculation gives
 \begin{equation}\label{for-A}
 A= e_3 \otimes \nabla u - \nabla u \otimes e_3  = 
 \left [ \begin{array}{ccc} 0 & 0 & - \partial_1 u \\
 0 & 0 & - \partial_2 u \\
 \partial_1 u & \partial_2 u & 0 
 \end{array}
 \right ]
 \end{equation} which implies
  \begin{equation}\label{for-B}
  \aoverbrace[L1R]{B} = \aoverbrace[L1R]{\frac{A^2}{2}} = -\frac12 \nabla u \otimes \nabla u.  
  \end{equation} However, we already know that  if $2<\beta<4$, then ${\rm sym} \nabla w = \aoverbrace[L1R]{B}$. Hence the proof of (iv) is complete. 
  
\medskip

 {In order to see (v), it is sufficient to repeat the same argument as in (iv) using the new assumptions. The weak convergence of ${\rm sym} \nabla w_h$ being replaced by strong convergence, we deduce the strong convergence of the adjusted sequence $\tilde w_h$. Meanwhile, a careful reexamination of the proof of   Corollary \ref{for-R-cor} alongside the observations made in  \cite[Proposition 4]{FJMhier} (see also \cite[Corollary 4.2]{CDS14}) implies that the family 
 \begin{equation*} 
\mathscr{Y}= \frac 1 {h^\beta} \Big (\frac 1h \int_0^h |\nabla y^h (\cdot, x_3)-   R^h|^2 dx_3\Big ) 
 \end{equation*}is equi-integrable over $\Omega$.  Now, since $\beta = 2
\delta$ for $\beta\ge 4$, the estimate 
$$
 \begin{aligned}
 \int_{B(x,\bar c h) \cap \Omega} |W_h|^2  \le Ch^2 \Big (  \frac 1h \int_{(B(x, \bar c h) \cap \Omega) \times (0,h)} |\nabla y^h -   R^h|^2 \Big ) 
 \end{aligned} 
$$ in combination with the equi-integrability of $\mathscr{Y}$ implies once again that 
$$
   \fint_{B(x,\bar c h) \cap \Omega} |w_h|^2  \xrightarrow{h\to0} 0.
 $$ The conclusion then follows similarly as in (iv). }
  \end{proof}

   \begin{lemma}\label{avg-to-0}
   Let $\Omega \subset \R^2$ be  a bounded Lipschitz  domain, $\bar c >0$ and assume that for a subsequence  $h_j\to 0$, $u_{h} \in W^{1,2}(\Omega)$  converge strongly to $0$ in $W^{1,2}$. Then,  up to a subsequence, and for q.e.~$x\in \overline \Omega$
   $$   
   \lim_{h\to 0} \fint_{B(x, \bar c h) \cap \Omega}  |u_h (y)| \, dy =0. 
   $$
   \end{lemma}
 \begin{proof}
Since $|u_h|$ also converge strongly to $|u|$, without loss of generality we can assume that $u_h \ge 0$ a.e. in $\Omega$. We first extend the $u_h$ to the whole $\R^2$ using the bounded linear extension operator $\mathcal E$  introduced in Appendix \ref{qeAppx}, formula \eqref{defEG}. {It directly follows from    Proposition \ref{propositionNice}  that
 $$
  \lim_{h\to 0} \fint_{B(x, \bar c h)} \mathcal{E}(u_h)(y) \, dy=0 \Rightarrow \lim_{h\to 0} \fint_{B(x, \bar c h) \cap \Omega}  u_h (y) \, dy =0.
 $$ 
 }
 
Therefore it is sufficient to prove the statement for  $\Omega = \R^2$.  In that case, we choose a suitable subsequence of $u_{h_j}$ which satisfies, after relabeling,  
$$
\int_{\R^2}  |Du_{h_j}|^2 \le \frac{1}{2^{3j}}. 
$$ and we define  
$$
B_j:= \Big \{x\in \R^2; \,\, \fint_{B(x, \bar ch_j)} u_{h_j}(y) \,dy > \frac 1{2^{j}} \Big \} 
\subset \Big \{x\in \R^2; \,\, \fint_{B(x, r)} u_{h_j}(y) \,dy > \frac 1{2^{j}} \,\, \mbox{for some} \,\, r>0 \Big \} .  
$$ Applying \cite[Lemma 4.8.1]{EG} to $u_{h_j}$, we obtain that 
$$
{\rm Cap}_{1,2}(B_j) \le C {2^{2j}} \int_{\R^2}  |Du_{h_j}|^2 \le \frac{C}{2^j}.
$$ Let
$$
E_k:=  \bigcup_{j=k}^\infty B_j. 
$$ Then 
$$
{\rm Cap}_{1,2} (E_k) \le  \frac{C}{2^{k-1}}
$$ and  for all $x\in \R^2 \setminus E_k$ we have for all $j\ge k$
$$
\fint_{B(x, \bar ch_j)} u_{h_j}(y)  \le \frac{1}{2^{3j}}.
$$ Thus by letting 
$$
E:= \bigcap_{k=1}^\infty E_k
$$ we have ${\rm Cap}_{1,2}(E)=0$ and 
$$
\forall x\in \R^2 \setminus E \quad \lim_{j\to \infty}\fint_{B(x, \bar ch_j)} u_{h_j}(y) =0. 
$$ The proof is complete. \end{proof}

\begin{cor}\label{strain} Let $\Omega$, $\Omega^h$, $K$, $\beta$ be as above. Assume that for a sequence as $h\to 0$, $y^h \in \mathcal{A}^h_K$, $E^h(y^h) \precsim h^\beta$,  $R^h$ is chosen as in \eqref{for-Rh}, and  $u$, $w$ are as in Theorem \ref{compactness}. Let 
$$
G^h:= \frac {(R^h)^T \nabla_h \tilde y^h -   {\rm Id}}{h^{\beta/2}}.
$$ Then up to a subsequence $\ds G^h  \rightharpoonup G$ weakly in $ L^2(\Omega^1, \R^{3\times 3})$ where 
\begin{equation}\label{Gcorner}
\aoverbrace[L1R]{G} = G_0  (x') + (x_3 -\frac 12) G_1(x'),
\end{equation} with
$$
G_1 = -\nabla^2 u
$$ and 
\begin{itemize}
\item[(i)] if $\beta=4$, then $\ds {\rm sym}\,G_0  = {\rm sym} \nabla w  + \frac 12 \nabla u \otimes \nabla u$, 
\item[(ii)]  if $\beta>4$, then  ${\rm sym} \,G_0  = {\rm sym} \nabla w $ 
  \end{itemize}
\end{cor}

\begin{proof}
The proof  is as in \cite[Lemma 2]{FJMhier}. We include it for the convenience of the reader. We define the difference quotient
$$
H^h(x', x_3,s):=  \frac {G^h(x', x_3+s) - G^h (x', x_3)}{s}
$$  Letting 
$$
\eta^h(x', x_3,s):= \frac{1}{h^{\beta/2-1}}  \Big ( \frac 1s \int_0^s  \frac 1h \partial_3 (\tilde y^h)'(x', x_3 + \sigma) d\sigma \Big ) 
$$
We first observe that for all $x_3, x_3 + s \in [0,1]$
$$
\begin{aligned}
\Big \| \eta^h(x', x_3,s)  &  - (A_{13}, A_{23})^T \Big \|_{L^2(\Omega)}   = 
\Big   \|\frac{1}{h^{\beta/2-1}}  \Big ( \frac 1s \int_0^{s}   \partial_3 (y^h)'(x',h (x_3 + \sigma) ) d\sigma  \Big )   - (A_{13}, A_{23})^T \Big \|_{L^2(\Omega)} 
 \\ & 
\le    \Big \|\frac{1}{h^{\beta/2-1}} \Big ( \frac 1s \int_0^s   \partial_3 (y^h)'(x', h (x_3 + \sigma) )  d\sigma  \Big ) 
  - (A^h_{13}, A^h_{23})^T \Big \|_{L^2(\Omega)}  + \|A^h - A\|_{L^2(\Omega)}
 \\ & 
 \le    \Big \|  \frac 1s \int_0^s   \frac{1}{h^{\beta/2-1}}  (\nabla y^h - {\rm Id}) (x', h (x_3 + \sigma) )   - A^h(x')    d\sigma   \Big \|_{L^2(\Omega)}  + \|A^h - A\|_{L^2(\Omega)} 
 \\ & 
 \le  \frac 1{\sqrt {hs}}    \Big \| \frac{1}{h^{\beta/2-1}}  (\nabla y^h - {\rm Id})  - A^h  \Big \|_{L^2(\Omega \times [hx_3 , h(x_3 + s)])}  + \|A^h - A\|_{L^2(\Omega)} 
  \\ & 
\le 
  \frac 1{\sqrt {hs}}    \Big \| \frac{1}{h^{\beta/2-1}}  (\nabla y^h - R^h)    \Big \|_{L^2(\Omega^h)}  + \|A^h - A\|_{L^2(\Omega)} \xrightarrow{h\to 0} 0,
\end{aligned}
$$ where we used \eqref{for-Rh} and the definition of $A^h$ in \eqref{for-Ah}.  As a consequence, and in view of \eqref{for-A}, 
$$
\eta^h (\cdot, x_3, s) \longrightarrow  -\nabla u \quad \mbox{in}\,\,  L^2(\Omega)  
$$ which implies 
$$
\nabla' \eta^h (\cdot, x_3 ,s ) \longrightarrow - \nabla^2 u  \quad \mbox{in}\,\, {\mathcal D}'(\Omega).
$$  On the other hand we have by straightforward calculation of the integral in the definition of $\eta^h$
$$
\eta^h (\cdot, x_3 ,s )  =  h^{-\beta/2} \frac{(\tilde y^h)' (\cdot,x_3 + s) - (\tilde y^h)' (\cdot,x_3)}{s},
$$ which yields
$$
\nabla' \eta^h (\cdot, x_3 ,s ) = \aoverbrace[L1R]{R^h H^h} (x', x_3, s) \rightharpoonup 
\frac {\aoverbrace[L1R]{G}(x', x_3+s) - \aoverbrace[L1R]{G}(x',x_3)}{s}    \quad \mbox{weakly in} \,\,L^2(\Omega), 
$$ since $R^h$ is uniformly bounded and  converges to ${\rm Id}$ in $L^2(\Omega)$. Comparing the two expressions of the limit for $\nabla' \eta^h$ we obtain
$$
\frac {\aoverbrace[L1R]{G} (\cdot, x_3+s) - \aoverbrace[L1R]{G} (\cdot,x_3)}{s}   = -\nabla^2 u
$$ Since $\nabla^2 u$ is independent of $x_3$, we obtain that $\aoverbrace[L1R]{G}$ is affine in $x_3$, and  can be written in the form given in \eqref{Gcorner}, and that  $G_1 = -\nabla^2 u$.  In order to determine ${\rm sym} \, G_0$, we let
$$
G^h_0 (x') :=\int_0^1 G^h(x', x_3) \, dx_3   
$$ Note that
$$
\begin{aligned}
G^h & =  \underbrace{\frac{1}{h^{\beta/2}} ({\nabla_h \tilde y^h   -   {\rm Id})} -  \frac{1}{h^{\beta/2}}  (R^h  - {\rm Id})}_{\widetilde G^h} + \underbrace{\frac{1}{h^{\beta/2}} (R^h - {\rm Id})^T   (\nabla_h \tilde y^h  -  R^h)}_{\overline G^h}.
\end{aligned}
$$  The last term $\overline G^h$ converges to 0 in $L^2(\Omega^1)$ in view of  \eqref{for-Rh}. 
Considering the upper left sub-matrix of the first two terms we have
$$
\begin{aligned}
  \int_0^1 {\rm sym} \aoverbrace[L1R]{\widetilde G^h} (\cdot , x_3) \, dx_3    
 & = \frac{1}{h^{\beta/2}}  \int_0^1 {\rm sym}  (\nabla' (\tilde y^h)'- {\rm Id}  )    (\cdot , x_3) \, dx_3    -     \frac{1}{h^{\beta/2}}{\rm sym} (\aoverbrace[L1R]{R^h} - {\rm Id})    
\\ & =  \frac{1}{h^{\beta/2}} {\rm sym}\nabla \Big (\frac 1h \int_0^h ((y^h)' - {\rm id})  (\cdot , x_3) \, dx_3 \Big )   -    \frac{1}{h^{\beta/2}}{\rm sym} (\aoverbrace[L1R]{R^h} - {\rm Id})     
\\ & = \frac{h^{\min\{\beta-2, \beta/2\}}}{h^{\beta/2}} {\rm sym} \nabla w_h  -  \frac{h^{\beta-2}}{h^{\beta/2}} \aoverbrace[L1R]{B^h}  \end{aligned}
$$ by \eqref{def-wh} and \eqref{symRh-2}.  As a consequence, and in view of  Theorem \ref{compactness}(iii), \eqref{symRh-2} and \eqref{for-B}
$$
{\rm sym}  \aoverbrace[L1R]{ G^h_0}\xrightarrow{\mbox{\tiny \it weakly in} \,\, L^2(\Omega)} \left \{ \begin{array}{ll} {\rm sym} \nabla w  + \ds \frac 12 \nabla u \otimes \nabla u     & \mbox{if} \,\, \beta=4 \\ {\rm sym} \nabla w   & \mbox{if} \,\, {\beta>4}  \end{array}  \right .
$$  But on the other hand the weak limit of ${\rm sym}  \aoverbrace[L1R]{ G^h_0}$ must be equal to $\ds \int_0^1 \aoverbrace[L1R]{G}(\cdot, x_3) \, dx_3$, which is equal to $G_0$ in view of \eqref{Gcorner}. The proof is complete. \end{proof}

\subsection{$\Gamma$-$\liminf$ estimates}
 \label{liminf}

       \begin{theorem}\label{liminf-thm}
 Let $\Omega$, $\Omega^h$, $K$ be as above, and assume that ${\rm Cap}_{1,2}(K)>0$ and $\beta>2$. Assume that for a sequence as $h\to 0$, $y^h \in \mathcal{A}^h_K$,  
 \begin{itemize}
 \item[(H)] $\ds u_h  \rightarrow u$, in $W^{1,2}(\Omega)$, and $\ds w_h \rightharpoonup  w$ weakly in $W^{1,2}(\Omega)$,
   \end{itemize}  where $u_h$ and $w_h$ are respectively defined  as in \eqref{def-uh} and \eqref{def-wh}.
 Then

 \begin{itemize}
 \item[(i)]  If $2<\beta<4$  
 $$
\liminf_h \frac{1}{h^\beta} J^h(y^h) \geq    \int_{\Omega} \frac{1}{24}Q_2(\nabla^2 u) - u f \; dx. 
 $$

 \item[(ii)]  If $\beta=4$

 $$
 \liminf_h \frac{1}{h^\beta} J^h(y^h) \geq   \int_{\Omega} \frac{1}{2}Q_2(  {\rm sym} \nabla w  + \frac 12 \nabla u \otimes \nabla u) + \frac{1}{24}Q_2(\nabla^2 u) - u f \; dx.
    $$ 
  
 \item[(iii)]  {If $\beta>4$
 
  $$
 \liminf_h \frac{1}{h^\beta} J^h(y^h) \geq   \int_{\Omega} \frac{1}{2}Q_2(  {\rm sym} \nabla w ) + \frac{1}{24}Q_2(\nabla^2 u) - u f \; dx.
    $$ }

     \end{itemize}
 \end{theorem}
 \begin{proof}
 
 We first note that under the hypothesis (H) and the scaling the linear term of the energies converge, i.e.\@
 $$
  -  \frac 1h \int_{\Om^h} {\bf f}^h \cdot y^h \; dx \xrightarrow{h\to 0}\int_{\Om} - u f \; dx
 $$
  Therefore, in case $\liminf_h \frac{1}{h^\beta} E^h(y^h)= +\infty$, the results follow trivially. Otherwise, for a suitable subsequence
and a constant $C>0$, we  obtain $E^h(y^h) \le Ch^\beta$, and hence the hypotheses of Theorem \ref{compactness} are satisfied. Note that the limits $u,w$ obtained in Theorem \ref{compactness} must be the same given under the convergence hypothesis (H), and so we obtain that $u\in H^2_{0,K}(\Omega)$. Applying Corollary \ref{strain}, the rest of the proof  follows exactly as in \cite[Corollary 2]{FJMhier}, and is left to the reader. \end{proof}

 \subsection{Recovery sequence and $\Gamma$-$\limsup$ estimate for $\beta>4$}   
  In this section we explain how to construct a recovery sequence. Namely, we want to prove the following statement.

 \begin{theorem}\label{limsup-thm}  {Let $\beta>4$. For every {$u \in H^2_{0,K}(\Omega)$}, $w\in H^1_{0,K}(\Omega)$, there exists a sequence $y^h \in  \mathcal{A}^h_K$ such that for the sequences $u_h, w_h$ defined as in \eqref{def-uh} and \eqref{def-wh}, 
  $(u_h, w_h) \longrightarrow (u,w)$ in $W^{1,2}(\Omega)$, and 
$$
\limsup_h \frac{1}{h^\beta} J^h(y^h) \leq  \int_{\Omega}\frac{1}{2}Q_2(  {\rm sym} \nabla w )  +  \frac{1}{24}Q_2(\nabla^2 u) - u f \; dx.
$$}
  \end{theorem}

  \begin{proof} { The proof starts  from the construction given in \cite[Section 6.2]{FJMhier}. More precisely we  take $\beta=2\alpha -2$ fix  given {$u\in H^2_{0,K}(\Omega)$, $w\in H^1_{0,K}(\Omega)$, } and $g\in L^2$ a given vector field. Then we define

  $$y^h(x',x_3)=\left(
  \begin{array}{c}
  x' \\
  h ( x_3 - \frac 12)
  \end{array}
  \right) + \left(
  \begin{array}{c}
  h^{\alpha-1} w \\
  h^{\alpha-2} u
  \end{array}
  \right) -h^{\alpha-1}(x_3 -
  \frac 12) \left(
  \begin{array}{c}
  \partial_1 u \\
 \partial_2 u \\
 0
  \end{array}
  \right)  +\frac{h^\alpha}{2}( x_3 - \frac 12)^2 g. $$
  
The required convergence of $(u_h, w_h)$  is established by a straightforward computation. Also, computations as in  \cite{FJMhier} imply that 
  $$\frac{1}{h^\beta}W(\nabla y^h)\xrightarrow{h \to 0} \frac{1}{2}Q_3({\rm sym} \nabla w + x_3B)  = \frac{1}{2}Q_3({\rm sym} \nabla w) +  \frac 12 x_3 ({\rm sym}\nabla w :  B) +  \frac{1}{2}Q_3(x_3B)  $$
   where  
  $$B:=-\nabla^2u+sym( g\otimes e_3).
  $$ As observed in \cite{FJMhier}, one can choose $g$ in such a way that 
  $$\int_{\Omega\times{(-\frac{1}{2}, \frac{1}{2})}}Q_3(x_3B)=\frac{1}{12}\int_{\Omega}Q_2(\nabla^2 u).$$ This choice yields
  $$
  \frac{1}{h^\beta} \int_{\Omega\times{(-\frac{1}{2}, \frac{1}{2})}} W(\nabla y^h) \xrightarrow{h \to 0}  \frac{1}{2} \int_\Omega  Q_3({\rm sym} \nabla w)  +  \frac{1}{24}\int_{\Omega}Q_2(\nabla^2 u).
  $$
  
  Now we see in the direct expression of $y^h$ and the fact that $u,w\in H^1_{0,K}(\Omega)$ implies that $y^h=0$ on $K$. This is not enough to have $y^h\in  \mathcal{A}^h_K$ because we need $y^h=0$ on $K_h$, a neighborhood of $K$. Therefore we need first to approximate $u,w$ by some $u_h, w_h$ which are   equal to $0$ on $K_h$. This can easily be done by use of a  sort of  diagonal argument. Indeed, {$u\in H^2_{0,K}(\Omega), w\in H^1_{0,K}(\Omega)$, and  thus from Lemma \ref{QeLemma0} we  know that there exist  sequences  $u_n, w_n \in C^\infty(\bar \Omega)$  such that $supp(u_n) \cup supp(w_n)\cap K=\emptyset$,  $u_n \to u$ in $W^{2,2}(\Omega)$, and $w_n \to w$ in $W^{1,2}(\Omega)$.} Since for $n$ fixed,  $K_h\subset \Omega \setminus ({\rm supp}(u_n) \cup {\rm supp}(w_n))$ for $h$ small,  we can then define a subsequence $n_h \to +\infty$ such that $K_h \subset \Omega \setminus  ({\rm supp}(u_{n_h}) \cup {\rm supp}(w_{n_h}))$. By this way if we change $y^h$ as above with $u_{n_h}, w_{n_h}$ instead of $u,w$, then we have now $y^h \in  \mathcal{A}^h_K$. Since the convergences of $u_{n_h}$ to $u$, and $w_{n_h}$ to $w$, hold respectively in the strong topologies of  $W^{2,2}$ and $W^{1,2}$,  it is easy to verify that we still have the same  limit 
   $$\frac{1}{h^\beta}W(\nabla y^h) \xrightarrow{h \to 0}  \frac{1}{2}Q_3({\rm sym} \nabla w + x_3B),$$
   which finishes the proof of the proposition.}   \end{proof}
  
  \begin{rem}\label{no-w}
 {Theorems \ref{liminf-thm} and \ref{limsup-thm} mean that for $\beta>4$, the functionals $h^{-\beta} E^h$, 
  $\Gamma$-converge to the functional
  $$
  I(u,w):=  \int_{\Omega}\frac{1}{2}Q_2(  {\rm sym} \nabla w )  +  \frac{1}{24}Q_2(\nabla^2 u) - u f \; dx,
$$  defined over $H^2_{0,K}\times H^1_{0,K}$, under the strong convergence of $(u_h, w_h)$ to $(u,w)$  in
$W^{2,2}\times W^{1,2}$.  Note that the strong convergence of $w_h$, and hence the inclusion of the weak limit $w$ in $H^1_{0,K}$ in Theorem \ref{compactness} is not established in its full generality for this regime. To bypass this problem, and in view of the fact that $u$ and $w$ are decoupled in the expression of $I(u,w)$,  one can forgo the mention of $w$ in the final result. Indeed, notice that both inequalities still work for $w=0$, i.e. 
$$
\liminf_h \frac{1}{h^\beta} J^h(y^h) \geq    \frac{1}{24}Q_2(\nabla^2 u) - u f \; dx,
$$ and the same recovery sequence written for $u\in H^2_{0,K}(\Omega)$ and $w=0$ satisfies 
$$
\limsup_h \frac{1}{h^\beta} J^h(y^h) \leq  \frac{1}{24}Q_2(\nabla^2 u) - u f \; dx.
$$ As stated in Theorem \ref{main11}, this establishes the functional 
$$
I(u) := \int_{\Omega} \frac{1}{24}Q_2(\nabla^2 u) - u f \; dx, 
$$ over $H^2_{0,K}$ as the $\Gamma$-limit of $h^{-\beta} E^h$  under the strong convergence of $u_h$ to $u$. 
}
  \end{rem}

  \subsection{Energy scalings and $\beta$-minimizing sequences}
   
 We first prove the following proposition regarding the fact that the scaling of the infimum energy is determined by the scaling of body forces in our setting.  Remember that vertical body forces ${\bf f}^h: \Om^h \to \R^3$ are assumed to be of the form ${\bf f}^h:= (0,0, h^\alpha \tilde f)$, for $\tilde f\in L^2(\Omega)$.

  \begin{prop}\label{energy-bnds}  Assume $\alpha > 2$, and let $\beta= 2\alpha-2$. There exists $h_0>0$ such that for a constant $C:= C(\Omega, K, \tilde f, h_0)>0$  and all $h<h_0$ 
  $$
-C h^{\beta} \le  \inf_{\mathcal {A}^{h}_K} J^h \le C h^{\alpha} \ll C h^{\beta}.
  $$   In particular, since $\inf_{\mathcal {A}^{h}_K} J^h \in \R$,  there exists a $\beta$-minimizing sequence. 
  \end{prop}
\begin{proof}
First note that letting $y = {\rm id}_{\Om^h}$ the identity map we have
$$
J^h({\rm id}_{\Om^h}) = - \frac 1h \int_{\Om^h} {\bf f}^h (x) \cdot x \, dx \le C h^\alpha \|\tilde f\|_{L^2(\Om)}
$$  yielding the upper bound on  $\inf_{\mathcal {A}^{h}_K} J^h $.

Let now $y\in {\mathcal {A}^{h}_K}$ be an arbitrary deformation, and consider the matrix field $\widetilde F$ obtained from Lemma \ref{mollified-grady}. Letting $\ds \widetilde F_0 := \fint_{\Omega} \widetilde F \, dx$.  We have by the Poincar\'e inequality
$$
\|\widetilde F-  \widetilde F_0\|^2_{L^2(\Omega)}  \le C \|\nabla \widetilde F\|^2_{L^2(\Omega)} \le \frac{C}{h^2} E^h(y),  
$$ which implies 
$$
\frac 1h \|\nabla y - \widetilde F_0\|^2_{L^2(\Omega^h)} \le  2 (\frac 1h \|\nabla y - \widetilde F\|^2_{L^2(\Omega^h)} 
+ \frac 1h\|F - \widetilde F_0\|^2_{L^2(\Omega^h)} ) \le  \frac{C}{h^2} E^h(y).
$$ Note that 
$$
|\widetilde F_0|^2 \le \fint _\Omega |\widetilde F|^2 \, dx \le  C (\fint _\Omega {\rm dist}^2 (\tilde F, SO(3)) \, dx + 1) \le C (\|\tilde d\|_{L^2(\Omega)}  + 1) \le C (E^h(y) + 1).
$$  Combining the last two estimates we obtain 
$$
\frac 1h \|\nabla y\|^2_{L^2(\Omega^h)}  \le C \Big ( \frac{1}{h^2} E^h (y) +1 \Big).
$$
Now, since $y\in \mathcal{A}^h_K$, we have $y_3 =0$ on $K_{h}$, i.e $y_3 \in {\bf A}^h_K$ and we can apply Theorem \ref{thin-poinc} to obtain
$$
\|y_3\|^2_{L^2(\Omega^h)}  \le  C \|\nabla y_3\|^2_{L^2(\Omega^h)}.
$$ Therefore
$$
\frac 1h \|y_3\|^2_{L^2(\Omega^h)}  \le \frac{C}{h} \|\nabla y_3\|^2_{L^2(\Omega^h)} \le C (\frac{1}{h^2} E^h(y) + 1), 
$$  and hence
$$
\begin{aligned}
J^h(y) & = E^h(y) - \frac 1h \int_{\Omega^h} {\bf f}^h \cdot y  \, dx   
\ge  E^h(y) - Ch^{\alpha} \|\tilde f\|_{L^2(\Omega)} \Big ( \frac 1{\sqrt h} \|y_3\|_{L^2(\Omega^h)} \Big )
\\ & \ge  E^h(y) - Ch^{\alpha} \|\tilde f\|_{L^2(\Omega)}  (\frac{1}{h} (E^h(y))^\frac 12 + 1) 
\\ & 
=  E^h(y)  - Ch^{\alpha-1}  \|\tilde f\|_{L^2(\Omega)} E^h(y)^\frac 12 - Ch^\alpha \|\tilde f\|_{L^2(\Omega)}. 
\end{aligned}
$$   Hence
\begin{equation}\label{JE}
J^h(y) \ge (E^h(y)^\frac 12 - Ch^{\alpha-1}\|\tilde f\|_{L^2(\Omega)}  )^2 - C^2h^{2\alpha-2}\|\tilde f\|^2_{L^2(\Omega)} - Ch^\alpha \|\tilde f\|_{L^2(\Omega)}. 
\end{equation}

A standard argument now yields for an adjusted constant $C$
$$
J^h(y) \ge -  Ch^{2\alpha-2}  \|\tilde f\|_{L^2(\Omega)} \ge -Ch^\beta,
$$  implying the required lower bound on  $\inf_{\mathcal {A}^{h}_K} J^h$. 
\end{proof}

\begin{cor}\label{beta-compact} 
Let $\alpha>2$ and $\beta= 2\alpha-2$.  For any sequence $y^h \in{\mathcal {A}^{h}_K}$,  
$$
J^h(y^h) \precsim h^\beta  \implies E^h(y^h) \precsim h^\beta.
$$ 
 In particular any $\beta$-minimizing sequence $y^h$ 
satisfies $E^h(y^h) \le Ch^\beta$.  
\end{cor}
\begin{proof}
Assuming $J^h(y^h) \precsim h^\beta$, from \eqref{JE} we obtain
$$
(E^h(y^h)^\frac 12 - Ch^{\beta/2} )^2 \le J^h(y^h) + Ch^\beta +Ch^\alpha \le Ch^\beta
$$ which implies$ E^h(y^h) \precsim h^\beta$.  Consider now a $\beta$-minimizing sequence $y^h$,  
$$
\displaystyle \limsup_{h \to 0}  \frac 1{h^\beta} \Big  (J^h (y^h) - \inf_{\mathcal {A}^{h}_K} J^h \Big ) =0.
 $$ We have therefore for $h$ small enough
$$
J^h (y^h) \le (\inf_{\mathcal {A}^{h}_K} J^h +1)  h^\beta  \le Ch^\beta,
$$ which implies the required estimate on $E^h(y^h)$. \end{proof}

\begin{theorem}
Let $\alpha>3$, $\beta = 2\alpha-2>4$,  and let $y^h \in \mathcal{A}^h_K$ be a $\beta$-minimizing sequence for $J^h$, and let $(u_h, w_h)$ be defined as in \eqref{def-uh}, \eqref{def-wh}. Then up to a subsequence, $(u_h, w_h) \longrightarrow (u,0)$ in $W^{1,2}(\Omega)$,  {$u\in H^2_{0,K}(\Omega)$}, and $u$ minimizes  the functional
$$
I(u) := \int_{\Omega} \frac{1}{24}Q_2(\nabla^2 u) - u f \; dx, 
$$ over {$H^2_{0,K}(\Omega)$.} Conversely, if  $u$ minimizes $I(u)$   over {$H^2_{0,K}(\Omega)$}, then there exists a 
$\beta$-minimizing sequence $y^h \in \mathcal{A}^h_K$ for $J^h$ such that $(u_h, w_h) \longrightarrow (u,0)$ in $W^{1,2}(\Omega)$.
 \end{theorem} 
 \begin{proof}
Combining Corollary \ref{beta-compact} and Theorem \ref{compactness}, we obtain that up to a subsequence $u_h \longrightarrow u$ in $W^{1,2}(\Omega)$, {$u\in H^2_{0,K}(\Omega)$}, and that $w_h \rightharpoonup w$ weakly in $W^{1,2}(\Omega, \R^2)$. A standard argument using the $\Gamma$-convergence results Theorem \ref{liminf-thm} and Theorem \ref{limsup-thm} (see also Remark \ref{no-w}) shows that  $u$ minimize the functional  
$
I(u)
$ over {$ H^2_{0,K}(\Omega)$}, and that 
$$
I(u) = \min_{u\in  {H^2_{0,K}(\Omega)}} I = \lim_{h\to 0} \frac1 {h^\beta}   J^h(y^h)= \lim_{h\to 0} \frac1 {h^\beta} \inf_{\mathcal{A}^h_K} J^h.
$$    On the other hand, Theorem \ref{liminf-thm} implies that the recovery sequence $y^h$ of Theorem \ref{limsup-thm} for $w=0$, which satisfies the required convergence criteria, satisfies
$$
\lim_{h\to 0} \frac 1 {h^\beta} J^h(y^h)   = I(u) = \lim_{h\to 0} \frac1 {h^\beta} \inf_{\mathcal{A}^h_K} J^h, 
$$ which implies that $y^h$ is $\beta$-minimizing. 
\medskip

{ In fact,  more can be shown with a finer reasoning. Following the same arguments as in  \cite[Section 7.2]{FJMhier}, the equi-integrability of the family $\mathcal{Y}$ and the strong convergence of ${\rm sym}\nabla w_h$  in $L^2$ can be established. As a consequence  the strong convergence of $w_h$ to  $w\in H^1_{0,K}$ follows immediately from Theorem \ref{compactness}-(v). Using this stronger convergence, and the full force of Theorems    \ref{liminf-thm} and  \ref{limsup-thm}, it can be shown that the pair $(u,w)$ must minimize the functional 
$$
\int_{\Omega} \frac{1}{2}Q_2(  {\rm sym} \nabla w )  + I(u)
$$ over $H^2_{0,K}\times H^1_{0,K}$. Since $u$ and $w$ are decoupled in the energy, we must have ${\rm sym} \nabla w=0$ for the minimizer $(u,w)$, and  hence the infinitesimal rigidity of displacements in $\R^2$ implies that $w(x)= Dx+ b$ is an affine map with $D\in so(2)$, skew-symmetric, and therefore $\det D\neq 0$ or $D=0$. Since $w\in H^1_{0,K}$ and ${\rm Cap}_{1,2}(K)>0$, we deduce that $w=0$, since otherwise $w$ can  vanish at only one single point. We conclude that $w_h \to 0$ in $W^{1,2}(\Omega)$.}  

 \end{proof}
  
  
   \section{The optimal biharmonic support problem}
      \label{recovery}


In this section we know focus on shape optimisation  Problem in \eqref{probmain2} { and we assume for simplicity that $\partial \Omega \subset K$. By this way the space $H^k_{0,K}(\Omega)$ reduces to the more classical $H^k_0(\Omega \setminus K)$}. Our aim is to prove existence and regularity of minimizers $K$ for that problem. 

 \subsection{Dual formulation}

Let $\Omega \subset \mathbb{R}^2$ be any open and simply connected domain. We denote by $\mathcal{K}(\Omega)$ the collection of all compact { connected}  subsets $K\subset \overline{\Omega}$ such that { $\partial \Omega \subset K$}. In particular, for all $K \in \mathcal{K}(\Omega)$, any connected component of $\Omega \setminus K$ is necessarily simply connected.

 Notice that for $u \in C^\infty_c(\Omega)$ it is very classical that
$$\int_{\Omega} |\nabla^2 u|^2 \;dx = \int_{\Omega} |\Delta u|^2 \;dx.$$

Indeed, a simple computation yields
\begin{eqnarray}
\int_{\Omega} |\Delta^2 u|^2 \;dx& =& \int_{\Omega} |\Delta u|^2 \;dx + 2\int_{\Omega}u_{1
2}^2 -u_{11}u_{22} \;dx, \notag \\
&=&  \int_{\Omega} |\Delta u|^2 \;dx + 2\int_{\Omega} {\rm div} F \;dx,
\end{eqnarray}
where 
$$F:=
\left(
\begin{array}{cc}
u_{12}u_2-u_{22}u_1 \\
u_{12}u_1 - u_{11}u_2
\end{array}
\right).
$$
Now 
$$\int_{\Omega} {\rm div} F \;dx=0$$
because $F\in C^\infty_c(\Omega)$. Therefore, the biharmonic equation can be equivalently solved either by minimizing $\int_{\Omega \setminus K} |\nabla^2 u|^2 \;dx$ or $\int_{\Omega \setminus K} |\Delta u|^2 \;dx$ in $H^2_0(\Omega \setminus K)$. 
We denote by $L^2_{sym}(\Omega)$ the space of symmetric matrix valued $L^2$ functions that we endow with the scalar product
$$\langle  A,M\rangle = \int_{\Omega} A:M \;dx.$$ 

We denote by ${\rm div} M$ the divergence of the Matrix valued function $M$, which consists in taking the divergence of each raw of $M$. A simple computation reveals that for   $u \in C^\infty_c(\Omega)$ and $M \in C^\infty_{sym}(\Omega)$ it holds
$$\int_{\Omega} \nabla^2 u : M \;dx = \int_{\Omega} u \; {\rm div} {\rm div}M \;dx,$$
which naturally extends to $u\in H^2_0(\Omega)$ and $M \in H^2_{sym}(\Omega)$.

We begin with an elementary Lemma. 
\begin{lemma}\label{lemmeS}
Let $f:E\times F\to \mathbb{R}$ be a real valued  function defined on two given sets $E,F$. Assume that $(u_0,v_0)$ is a saddle point, i.e. satisfies 
$$f(u_0,v_0)=\max_{v} f(u_0,v)=\min_u f(u,v_0).$$
Then
$$\inf_u \sup_v f(u,v)=\sup_v \inf_u f(u,v)=f(u_0,v_0).$$
\end{lemma}
 \begin{proof} For any $(u,w)$ fixed it is clear that  $f(u,w)\leq \sup_{v}  f(u,v)$ thus taking first
 the inf in $u$ and then sup in $w$ yields
  $$ \sup_w \inf_u f(u,w)\leq \inf_{u}\sup_{v}  f(u,v).$$
  For the reverse inequality we use that $(u_0,v_0)$ is a saddle point. Indeed,
  $$\inf_u \sup_v f(u,v)\leq \max_v f(u_0,v)= \min_{u} f(u,v_0)\leq \sup_v \inf_u f(u,v),$$
  which concludes the proof.
 \end{proof}

\begin{prop} \label{saddle2}Let $F:H^2_0(\Omega \setminus K)\times L^2_{sym}(\Omega)\to \mathbb{R}$ be defined by 
$$F(u,M):=2 \int_{\Omega} \nabla^2 u : M \; dx-\int_{\Omega} |M|^2 \; dx -2 \int_{\Omega} u f \;dx.$$
Then $(u_K, \nabla^2u_K)$ is a sadle point for $F$. In particular,
$$- \int_{\Omega}|\nabla^2u_K|^2\;dx=\sup_{M\in L^2_{sym}(\Omega)} \inf_{u\in H^2_0(\Omega \setminus K)} F(u,M).$$
\end{prop}
\begin{proof} From the elementary equality $|A+B|^2=|A|^2+|B|^2+2 A:B$ valid for any matrix $A,B$ we infer that 
$$|A|^2\geq 2 A:B - |B|^2,$$
and in particular for any $M\in L^2_{sym}(\Omega)$, we always have 
$$ F(u_K, \nabla^2u_K)\geq F(u_K,M),$$
and since the equality occurs for $M=\nabla^2u_{K}$ we can affirm 
\begin{eqnarray}
 F(u_K, \nabla^2u_K)=\max_M F(u_K,M). \label{saddle1}
 \end{eqnarray}

On the other hand, by Lax-Milgram theory we know that 
  $u_K$ is the unique minimizer in $H^2_0(\Omega\setminus K)$ for the functional $J$ defined by
$$J(u):= \int_{\Omega} |\nabla^2 u| dx -2 \int_{\Omega} u f \;dx,$$
and the weak formulation for this problem says that 
$$ \int_{\Omega} \nabla^2 u:\nabla^2u_K \; dx = \int_{\Omega} u f \;dx \quad \quad \forall u \in H^2_0(\Omega \setminus K).$$

In turn, for $\nabla^2u_K$ fixed and an arbitrary $u \in H^2_0(\Omega \setminus K)$, we observe that  the expression of $F(u,\nabla^2u_{K})$ reduces to 
$$F(u,\nabla^2u_{K})=-\int_{\Omega} |\nabla^2u_K|^2 \;dx,$$
which is constant in the $u$ variable. We deduce a fortiori that 
$$F(u_K,\nabla^2u_{K})= \min_u F(u,\nabla^2u_{K}),$$
which together with \eqref{saddle1},  proves that $(u_K,\nabla^2u_K)$ is a saddle point for $F(u,M)$. We conclude by applying Lemma \ref{lemmeS}.
\end{proof}

In the following proposition we will use the notation $LD(A)$  for a function $v\in L^2(A, \R^2)$ satisfying $e(v)\in L^2(A)$, where $e(v)=(\nabla v+\nabla v^T)/2$ is the symmetrized gradient of $v$.

\begin{prop} \label{dual}  Let $\varphi \in H^1_0(\Omega)$  be the unique solution of
$$
\left\{
\begin{array}{c}
-\Delta \varphi=f \text{ in }\Omega\\
\varphi \in H^1_0(\Omega)
\end{array}
\right.
$$
and let 
$$G=\varphi{\rm Id}=\left(
\begin{array}{cc}
\varphi & 0\\
0 & \varphi
\end{array}
\right).$$

Let $K$ be a minimizer for Problem~\eqref{probmain2}. Then there exists $v_K\in LD(B(x,r)\setminus K)$ such that  $(v_K, K)$ is a minimizer for the problem

\begin{equation} \label{3.2}
  \min_{(v, K) \in \mathcal{A}}  \displaystyle\int_{\Omega \setminus K} |e(v)- G|^{2} dx +   \mathcal{H}^{1}(K)
\end{equation}
where 
\begin{align*}
\mathcal{A}:=\{ (v,  K):  K \in \mathcal{K}(\Om)\, \text{ and } 
v \in LD(B(x,r)\setminus K ) \}.
\end{align*}
\end{prop}
\begin{proof} For a given $K\subset \overline{\Omega}$, we can apply Proposition \ref{saddle2}  to write
\begin{eqnarray}
-\int_{\Omega}|\nabla^2u_K|^2 \;dx&=&\sup_{M \in L^2_{sym}}  \inf_{u \in H^2_0(\Omega\setminus K)}   \; 2 \int_{\Omega} \nabla^2 u : M \; dx-\int_{\Omega} |M|^2 \; dx -2 \int_{\Omega} u f \;dx \notag \\
&=& \sup_{M \in L^2_{sym}}  \inf_{u \in H^2_0(\Omega\setminus K)}   \; 2  \;_{H^2_0}\langle   u , \; {\rm div }\,{\rm div }\,M \rangle_{(H^2_0)'}-\int_{\Omega} |M|^2 \; dx -2 \int_{\Omega} u f \;dx .\notag \\
&=& \sup_{M \in L^2_{sym}}  \inf_{u \in H^2_0(\Omega\setminus K)}   \; 2  \;_{H^2_0}\langle   u , \; {\rm div }\,{\rm div }\,M -f\rangle_{(H^{2}_0)'}-\int_{\Omega} |M|^2 \; dx   .\notag 
\end{eqnarray}
The infimum in the $u$ variable in the above  is equal to $-\infty$, excepted when ${\rm div }\,{\rm div }\,M=f$ in $\mathcal{D}'(\Omega \setminus K)$. This leads to the following dual equality
$$
\int_{\Omega} |\nabla^2u_{K}|^2 \;dx= \min \big\{   \int_{\Omega} |M|^2 \; dx    \quad \text{ s.t. }  M\in L^2_{sym}(\Omega) \text{ and } {\rm div }\,{\rm div }\,M=f \text{ in } \mathcal{D}'(\Omega \setminus K) \big \}.
$$
Moreover, the minimium is   achieved for $M=\nabla^2u_{K}$. In other words

$$\min_{K}\int_{\Omega} |\nabla^2u_{K}|^2 \;dx+\mathcal{H}^1(K)=\min_{K}\min_{M}\int_{\Omega} |M|^2 \;dx+\mathcal{H}^1(K),$$
where the minimum in $M$ is over all $M\in L^2_{sym}(\Omega)$ satisfying  ${\rm div }\,{\rm div }\,M=f \text{ in } \mathcal{D}'(\Omega \setminus K)$.

Now let $v \in H^1_0(\Omega)$  be the unique solution of
$$
\left\{
\begin{array}{c}
-\Delta v=f \text{ in }\Omega\\
v \in H^1_0(\Omega)
\end{array}
\right.
$$
and let 
$$G=v{\rm Id}=\left(
\begin{array}{cc}
v & 0\\
0 & v 
\end{array}
\right).$$

Then $G$ is symetric and $ {\rm div }\,{\rm div }\,G=\Delta v = -f $ in $\mathcal{D}'(\Omega)$. It follows that for all $M$ satisfying $ {\rm div }\,{\rm div }\,M=f \text{ in }\mathcal{D}'(\Omega \setminus K),$
$$ {\rm div }\,{\rm div }\,(M+G)=0 \text{ in }\mathcal{D}'(\Omega \setminus K).$$
In particular, since $\Omega \setminus K$ is simply connected, there exists $u$ such that 
$${\rm div }\,(M+G) = \nabla^\perp u.$$
We deduce that 
$${\rm curl ({\rm div}(M+G)^\perp)}=0$$
or differently
$${\rm curl curl} (Com(M+G))=0.$$
The classical Saint-Venant compatibility condition yields the existence of $u$ such that $Com(M+G)=e(u)$ in the simply connected domain $\Omega\setminus K$ and we infer that 
$$M=Com(e(u))-G,$$
and finally since $Com(G)=G$,
$$\int_{\Omega \setminus K}|M|^2\;dx=\int_{\Omega \setminus K} |e(u)-G|^2 \;dx.$$

We conclude that 
$$\min_{M}\int_{\Omega \setminus K}|M|^2\;dx=\min_{u}\int_{\Omega \setminus K} |e(u)-G|^2 \;dx$$
 which ends the proof of the proposition.
\end{proof}

\begin{rem} \label{bounded}Notice that  $G$ is bounded and we actually have 
$$\|G\|_{L^\infty(\Omega)}\leq C \|f\|_p.$$
Indeed,  by elliptic regularity $G \in W^{2,p}(\Omega)$ and then since $p>2$ we know that $u\in C^{1,\alpha}(\overline{\Omega})$. 
\end{rem}


\subsection{Existence}

We start by proving the existence of a minimizer, as a simple consequence of Sverak's \cite{sverak} result. 

\begin{prop} \label{existence}  Let $\Om$ be an open and bounded set in $\mathbb{R}^{2}$ and  let $f \in L^{\infty}(\Om)$. Let $({K}_{n})_{n}$ be a sequence of closed connected  subset of $\overline{\Om}$, converging to a closed connected set $K \subset \overline{\Om}$ with respect to the Hausdorff distance. Then
\[
u_{K_{n}} \xrightarrow{n \to +\infty} u_{K}\,\ \text{strongly in}\,\ H^{2}(\Om).
\]
\end{prop}
\begin{proof}  We proceed as a standard $\Gamma$-convergence argument. We start by noticing that for all $n$, using that $u_{K_n}$ is the solution of $\Delta^2 u_{K_n}=f$ in $H^2_0(\Omega\setminus K_n)$,
$$\int_{\Omega} |\nabla^2 u_{K_n}|^2 \;dx = \int_{\Omega} f u_{K_n} \leq \|f\|_2^2 \|u_{K_n}\|_{L^2}\leq C\|\nabla^2u\|_2,$$
where we have used the Poincar\'e inequality in $H^2_0(\Omega)$. This leads to 
$$\int_{\Omega} |\nabla^2 u_{K_n}|^2 \;dx \leq \|f\|_2, $$
so that $u_{K_n}$ is uniformly bounded in $H^2_0(\Omega)$. We can therefore extract a subsequence (not relabeled) that converges weakly in $H^2(\Omega)$, strongly in $H^1_0(\Omega)$, and uniformly in $\Omega$ to some function $v \in H_0(\Omega)$. By uniform convergence we also know that $v=0$ on the set $K$, the Hausdorff limit of $K_n$.

Since $\nabla u_{K_n} \in H^1_0(\Omega \setminus K_n)$ and  converges weakly in $H^1_0(\Omega)$ to $\nabla v$, and since the sets $K_n$ are all compact and connected in $\overline{\Omega}$, it follows from \cite{sverak} that $\nabla v \in H^1_0(\Omega)$. In other words $v \in H^2_0(\Omega \setminus K)$.

Now let $\varphi \in C^\infty_c(\Omega \setminus K)$. By Hausdorff convergence of $K_n$ we deduce that the support of $\varphi$ is contained in $\Omega \setminus K_n$ for all $n$ large enough. Therefore, we can apply the weak formulation of the problem satisfied by $u_{K_n}$ which yields 
$$\int_{\Omega} \nabla^2 u_{K_n} : \nabla^2 \varphi \;dx =\int_{\Omega} f \varphi \; dx.$$
Passing to the limit we obtain that $u_{K_n}$  is the unique solution of $\Delta^2 u =f$ in $H^2_0(\Omega \setminus K)$. In other words $u=u_{K}$.

By the week convergence of the Hessians we already have 
$$\int_{\Omega} |\nabla^2u| \;dx \leq \liminf \int_{\Omega} |\nabla^2 u_{K_n}|^2 \;dx.$$

Conversely, by definition of $H^2_0$ there exists a sequence of $v_n \in C^\infty_c(\Omega \setminus K)$ such that $v_n\to u$ strongly in $H^2(\Omega)$. Let us first fix a $k_0\geq 0$. Using the Hausdorff convergence of $K_n$ to $K$ we know that $v_{k_0} \in H^2_0(\Omega \setminus K_n)$ for $n$ large. Since $u_{K_n}$ is a minimizer in this class we obtain
$$
\int_{\Omega} |\nabla^2 u_{K_n}|^2 \; dx - 2\int_{\Omega} u_{K_n }f \;dx \leq \int_{\Omega} |\nabla^2 v_{k_0}|^2 \; dx - 2\int_{\Omega} v_{k_0 }f \;dx.
$$
passing to the limsup in $n$ we arrive at 
$$
\limsup_{n\to +\infty}\int_{\Omega} |\nabla^2 u_{K_n}|^2 \; dx \leq   \int_{\Omega} |\nabla^2 v_{k_0}|^2 \; dx - 2\int_{\Omega} v_{k_0 }f \;dx
+ 2\int_{\Omega} u f \;dx.$$
Letting now $k_0 \to +\infty$ yields
$$
\limsup_{n\to +\infty}\int_{\Omega} |\nabla^2 u_{K_n}|^2 \; dx \leq   \int_{\Omega} |\nabla^2 u|^2 \; dx ,$$
which together with the liminf inequality proves the strong convergence of $u_{K_n}$ in $H^2(\Omega)$ to $u=u_{K}$ (the convergence of the full sequence follows from the uniqueness of the limit). 
\end{proof}

 \begin{prop} \label{prop: 2.15}
Problem~\eqref{probmain2} admits a minimizer { in the class $\mathcal{K}(\Omega)$}.
\end{prop}
\begin{proof} Let $({K}_{n})_{n}$ be a minimizing sequence for Problem~\ref{probmain2} { in  $\mathcal{K}(\Omega)$}.  Using Blaschke's theorem (see \cite[Theorem 6.1]{APD}), we can find a compact connected proper subset $K$ of $\overline{\Om}$ such that up to a subsequence, still denoted by the same index, ${K}_{n}$ converges to $K$ with respect to the Hausdorff distance. Then, { of course $\partial \Omega\subset K$ and } by Proposition~\ref{existence}, $u_{{K_{n}}}$ converges to $u_{K}$ strongly in $H^2_{0}(\Om)$. Finally, thanks to the lower semicontinuity of $\mathcal{H}^{1}$ with respect to the  Hausdorff convergence { of connected sets}, we deduce that $K$ is a minimizer of Problem~\eqref{probmain2}.
\end{proof}


\subsection{Ahlfors regularity}
We recall that a set $K \subset \mathbb{R}^{2}$ is said to be Ahlfors regular  of dimension 1, if there exist some constants $c>0$, $r_{0}>0$ and $C>0$ such that for every $r \in(0, r_{0})$ and for every $x \in K$ the following holds
\[
cr\leq \mathcal{H}^{1}(K \cap B_{r}(x)) \leq Cr.\numberthis \label{5.1}
\]

Note that for a closed connected nonempty set $K$ the lower bound in (\ref{5.1}) is trivial: indeed, for all $x \in K$ and for all $r\in (0,\diam(K)/2)$ we have $K \cap \partial B_{r}(x) \not = \emptyset$, and then
\begin{equation}
\mathcal{H}^{1}(K \cap B_{r}(x)) \geq r.  \label{A}
\end{equation}

\begin{theorem} \label{thm: 5.3} Let $\Om \subset \mathbb{R}^{2}$ be a bounded Lipschitz domain. Let  $(v,K)$ be a solution of Problem~\ref{3.2} with $\diam(K)>0$. Then $K$   is Ahlfors regular. More precisely, there exists $C>0$ and $r_0>0$ such that for all $x\in K$ and $0\leq r \leq r_0$ we have 
$$\int_{B_r(x)\cap \Omega}|e(u)-G|^2 \,dx + \mathcal{H}^1(K \cap B_r(x))\leq 2\pi r +C\|f\|_pr^2.$$
 \end{theorem}

\begin{proof} It easily follows from the fact that $\Omega$ is a Lipschitz domain, that there exists $r_0>0$ such that $\partial B(x,r) \cap \overline{\Omega}$ is connected for all $x\in \Omega$ and $r\leq r_0$. Then to prove the Theorem, it suffice to compare $(u,K)$ with  the admissible competitor $(w,K_r)$ defined by 
\begin{equation}
K_{r}=(K \backslash B_{r}(x)) \cup (\partial B_{r}(x)\cap \overline{\Omega}),\label{3.6}
\end{equation}
and $w=u{\bf 1}_{\Omega \setminus B_r(x)}$.  By this way we obtain
$$\int_{B_r(x)\cap \Omega}|e(u)-G|^2 \,dx + \mathcal{H}^1(K \cap B_r(x))\leq \int_{B_r(x)\cap \Omega}|G|^2 \,dx +2\pi r,$$
and finally using also Remark \ref{bounded} we deduce that 
$$\int_{B_r(x)\cap \Omega}|e(u)-G|^2 \,dx + \mathcal{H}^1(K \cap B_r(x))\leq C\|f\|_p r^2+2\pi r,$$
which ends the proof.
\end{proof}


 \subsection{$C^1$ regularity}

The $C^1$ regularity of minimizers will follow from the following observation.

 \begin{prop} If $(u,K)$ is a minimizer of Problem \eqref{3.2}, then it is an almost-minimizer of the Grifith Energy. In other words there exists $C>0$ and $r_0>0$   such that for every competitor $(v,K')$ in the ball $B_r(x)$ and for all $r\leq r_0$ we have
  $$\int_{B_r}|e(u)|^2\;dx + \mathcal{H}^1(K \cap B_r(x)) \leq \int_{B_r}|e(v)|^2\;dx + \mathcal{H}^1(K' \cap B_r(x)) + Cr^{1+\frac{1}{2}}$$
 \end{prop}
 
 \begin{proof}   Let    $(v,K')$ be a competitor for $(u,K)$ in the ball $B_r(x)$. Then since $u$ is a minimizer for Problem \eqref{3.2} and $(v,K')$ coincides with $(u,K)$ outside $B_r$ we have
 \begin{eqnarray}
 \int_{B_r} |e(u)-G|^2 \,dx +\mathcal{H}^1(K \cap B_r) \leq  \int_{B_r} |e(v)-G|^2 \,dx +\mathcal{H}^1(K' \cap B_r) \notag
 \end{eqnarray}
 which implies 
  \begin{eqnarray}
 \int_{B_r} |e(u)|^2 \;dx +\mathcal{H}^1(K \cap B_r) \leq  \int_{B_r} |e(v)|^2   - 2\int_{B_r} e(v):G\,dx+ 2\int_{B_r} e(u):G\,dx+\mathcal{H}^1(K' \cap B_r). \notag
 \end{eqnarray}
Now by elliptic regularity (see Remark \ref{bounded}) we have 
$$\int_{B_r}|G |^2 \;dx \leq Cr^2,$$
thus after Cauchy-Schwarz we arrive at
  \begin{eqnarray}
 \int_{B_r} |e(u)|^2 \;dx +\mathcal{H}^1(K \cap B_r) \leq  \int_{B_r} |e(v)|^2 +   Cr \left(\int_{B_r} |e(v)|^2 \right)^{\frac{1}{2}} + Cr \left(\int_{B_r} |e(u)|^2 \right)^{\frac{1}{2}}  +\mathcal{H}^1(K' \cap B_r). \notag
   \end{eqnarray}

Now by the proof of Ahlfors-Regularity we   know that 
$$\int_{B_r} |e(u)-G|^2 \;dx+\mathcal{H}^1(K \cap B_r) \leq C_A r,$$
which implies in particular that 
$$\int_{B_r} |e(u)|^2 \;dx \leq 2\int_{B_r} |e(u)-G|^2 \;dx +2 \int_{B_r} |G|^2 \;dx \leq 2C_A r+Cr^2\leq C'r,$$
provided that $r_0\leq 1$ (that we can assume). Returning to the inequality obtained before, we infer that
$$ \int_{B_r} |e(u)|^2 \;dx +\mathcal{H}^1(K \cap B_r) \leq  \int_{B_r} |e(v)|^2 +    Cr \left(\int_{B_r} |e(v)|^2 \right)^{\frac{1}{2}} + Cr^{1+\frac{1}{2}} +\mathcal{H}^1(K' \cap B_r).$$

We now divide the argument in two alternatives. Either 
$$\int_{B_r} |e(v)|^2 \geq (C'+C_A) r$$ 
and then 
 $$\int_{B_r} |e(u)|^2 \;dx+\mathcal{H}^1(K \cap B_r) \leq (C'+C_A) r \leq \int_{B_r} |e(v)|^2 $$
 Or $\int_{B_r} |e(v)|^2 \leq (C'+C_A) r$ but then 
$$
 \int_{B_r} |e(u)|^2 \;dx +\mathcal{H}^1(K \cap B_r) \leq  \int_{B_r} |e(v)|^2 +   Cr^{1+\frac{1}{2}} +\mathcal{H}^1(K' \cap B_r). \notag
$$
 In both cases we always have 
 $$
 \int_{B_r} |e(u)|^2 \;dx +\mathcal{H}^1(K \cap B_r) \leq  \int_{B_r} |e(v)|^2   +\mathcal{H}^1(K' \cap B_r) + Cr^{1+\frac{1}{2}}. \notag
$$
 which achieves the proof of the Proposition.
 \end{proof}

Therefore, by the results contained in  \cite{lem-lab0} and \cite{lem-lab}, (see also  \cite{BIL}), we get the following theorem.

\begin{theorem} Let $K\in \mathcal{K}(\Omega)$ be a solution for Problem \eqref{3.2}. Then $K \cap \Omega$ is $C^{1,\alpha}$ outside a set of Hausdorff dimension strictly less than one. 
\end{theorem}

 
  \appendix

 
     \section{Quasi-everywhere traces of $W^{k,2}$ functions}\label{qeAppx}

 {

In the following statement we say that a property holds $(m,p)$-q.e. if it holds true outside a set of zero  
capacity of order $m\geq 1$ and integrability exponent $p$. We refer to \cite{ah} for the definition of  capacity ${\rm Cap}_{m,p}$.

Let $K\subset \R^N$  be a closed set and $u \in W^{k,p}(\R^N)$. For $1<p<+\infty$ there exists a nice characterization of the space $W^{k,p}_0(\R^N\setminus K)$ in terms of traces of $u$ on $K$ as follows. If $\alpha$ is a multiindex such that $|\alpha| \leq k$ we say that $\partial^\alpha u|_K=0$   if 
$$ \lim_{r\to 0} \left( \fint_{B_r(x_0)}|\partial^\alpha u| \; dx\right) = 0 \; \text{ for  $(k-|\alpha|,p)$-q.e. } x_0 \text{ on }K.$$
Then \cite[Theorem 9.1.3.]{ah} says that 
\begin{eqnarray}
W^{k,p}_0(\R^N \setminus K) =\left\{u \in W^{k,p}(\R^N) \text{ s.t. } \partial^{\alpha}u|_{K}=0  \text{ for all  } 0\leq \alpha \leq k-1\right\}, \label{h100}
\end{eqnarray}
 where here  $W^{k,p}_0(\R^N \setminus K)$ denotes the closure of $C^\infty_c(\R^N\setminus K)$ with respect to the $W^{k,p}(\R^N)$ norm.
The following lemma is a generalization of the above fact, adapted to a domain $\Omega$ instead of $\R^N$.   The proof uses  an  extension operator for $\Omega \subset \R^N$, i.e. a linear mapping $\mathcal{E}:W^{k,p}(\Omega)\to W^{k,p}(\mathbb{R}^N)$ such that $\mathcal{E}(u)=u$ in $\Omega$ and 
\begin{eqnarray*}
\|\mathcal{E}(u)\|_{W^{k,p}(\mathbb{R}^2)}\leq C\|u\|_{W^{k,p}(\Omega)}.  
\end{eqnarray*}
 The construction of such operator for a Lipschitz domain $\Omega$, is classical (see for instance \cite[Chapter 5]{adamsS}). By Lipschitz domain we mean a bounded open set whose boundary is locally the graph of a Lipschitz function at every point of the boundary.

 The strategy for a domain $\Omega$ is then to apply the characterization in \eqref{h100} to the extended function $\mathcal{E}(u)$ defined on $\R^N$, but it is not straightforward that the trace of $u$, which is a function defined only on $\Omega$, coincides with the trace of $\mathcal{E}(u)$ on the boundary $\partial \Omega$. For a Lipschitz domain this happens to be true thanks to a result in the book \cite{jonsson}.  Here  is then a general statement that we can prove with this strategy, that for simplicity we write only in the particular case $p=2$.

\begin{lemma}\label{QeLemma0}
Let $\Omega \subset \R^N$ be a bounded Lipschitz domain and $K\subset \overline{\Omega}$ be a closed set. We consider the subspace of $W^{k,2}(\Omega)$ defined by 
$$H^k_{0,K}(\Omega):=\left\{u \in W^{k,2}(\Omega) \text{ such that } \partial^{\alpha}u|_{K}=0  \text{ for all multiindices $\alpha$ such that } 0\leq \alpha \leq k-1\right\},$$
where  by $\partial^{\alpha}u|_K =0$ we mean that 
$$ \lim_{r\to 0} \left( \fint_{B_r(x_0) \cap \Omega}|\partial^\alpha u| \; dx\right) = 0 \; \text{ for  $(k-|\alpha|,2)$-q.e. } x_0 \text{ on }K.$$
Then  
\begin{enumerate}
\item $H^k_{0,K}(\Omega)=\big\{u \in W^{k,2}(\Omega) \text{ such that } \mathcal{E}(u) \in W^{k,2}_0(\R^N \setminus K)\big\}$.
\item $H^k_{0,K}(\Omega)\subset W^{k,2}(\Omega)$ is closed for the strong topology of $W^{k,2}(\Omega)$.
\item $H^k_{0,K}(\Omega)\subset W^{k,2}(\Omega)$ is closed for the weak topology of $W^{k,2}(\Omega)$.
\item If $u \in  H^k_{0,K}(\Omega)$ then there exists a sequence $\varphi_n \in C^{\infty}_c(\mathbb{R}^2)$ such that $supp(\varphi_n)\cap K= \emptyset$ for all $n$ and $\varphi_n \to u$ strongly in $W^{k,2}(\Omega)$.
\end{enumerate}

\end{lemma}

%

\begin{proof}  $\bullet$ {Proof of $(1).$}  We begin with the proof of $(1)$ that will actually be   the key ingredient for all the other points.  We already know by   \cite[Theorem 9.1.3.]{ah}   that 
\begin{eqnarray}
W^{k,2}_0(\R^N \setminus K) =\left\{u \in W^{k,2}(\R^N) \text{ s.t. } \partial^{\alpha}u|_{K}=0  \text{ for all  } 0\leq \alpha \leq k-1\right\}. \label{h10}
\end{eqnarray}
Thus to conclude, we only need to prove that for $(k-|\alpha|,2)$-q.e. $x_0 \in \partial \Omega$,
$$ \lim_{r\to 0} \left( \fint_{B_r(x_0) \cap \Omega}|\partial^\alpha u| \; dx\right) =   \lim_{r\to 0} \left( \fint_{B_r(x_0)}  |\partial^\alpha \mathcal{E}(u)| \; dx\right).$$
This fact actually follows from    \cite[Proposition 2 page 207]{jonsson}.   To say just a few words, the proof in \cite{jonsson} uses  that $\partial^\alpha(\mathcal{E}(u))$, as a global Sobolev function on the full space $\R^N$, can be written as a convolution with a potential, namely $v:=\partial^\alpha(\mathcal{E}(u))=G_\alpha * f$ where $f\in L^2$. Then  by the continuity behavior of Bessel potentials proved by Meyers in  \cite[Theorem 3.2.]{meyers} we know that this type of function has a sort of pointwise continuity property, in the sense that 
$$v(x_0)=\lim_{\substack{x \to x_0 \\ x \not \in E}} v(x)$$ 
where $E$ has a controlled capacity of the form ${\rm Cap}_{m,p}(E\cap B(x_0,r))=o(r^{N-1})$. From this we deduce that the limit of averages intersected with $\Omega$ or without intersection with $\Omega$ must coincide  $(k-|\alpha|,2)$-q.e. We refer to   \cite[Proposition 2 page 207]{jonsson} for more details. This achieves the proof of (1), and we can notice that it does not depend on the choice of the extension operator. Notice also that the argument used in Proposition \ref{propositionNice} could  give an independent proof for the special case $k=1$.

\medskip

$\bullet$ {Proof of $(2).$}  Let $u_n$ be a sequence in $H^k_{0,K}(\Omega)$ such that $u_n\to u$ in  $W^{k,2}(\Omega)$. Since the extension operator $\mathcal{E}$ is continuous on $W^{k,2}$ it follows that $\mathcal{E}(u_n)\to \mathcal{E}(u)$ in $W^{k,2}(\R^2)$ and by use of (1) we know  that $\mathcal{E}(u_n)\in W^{k,2}_0(\R^2\setminus K)$, for all $n$. But now by definition, $W^{k,2}_0(\R^N\setminus K)$ is a closed subspace of $W^{k,2}(\R^N\setminus K)$ and since $\mathcal{E}(u_n)\to \mathcal{E}(u)$ we obtain $\mathcal{E}(u)\in W^{k,2}_0(\R^2 \setminus K)$. By applying (1) again we deduce that $u\in H^k_{0,K}(\Omega)$.\\

 $\bullet$ {Proof of $(3).$} Let $u_n \in  H^k_{0,K}(\Omega)$ be a sequence that converges weakly in $W^{k,2}(\Omega)$ to some limit $u$, in other words $\langle u_n , \varphi \rangle\to \langle u,\varphi\rangle$ for all $\varphi \in W^{k,2}(\Omega)$, where the brackets means the scalar product in $W^{k,2}(\Omega)$. Since $H^k_{0,K}(\Omega)$ is a closed subset of $W^{k,2}(\Omega)$ (for the strong topology), it follows that  $H^k_{0,K}(\Omega)$ is itself a Hilbert space endowed with the same norm and scalar product of $W^{k,2}(\Omega)$. We can therefore extract a subsequence that converges for the weak topology  in $H^k_{0,K}(\Omega)$ to some limit function $h$ that must still belong to $H^k_{0,K}(\Omega)$. This says in particular that  $\langle u_n , \varphi \rangle\to \langle h,\varphi\rangle$ for all $\varphi \in H^k_{0,K}(\Omega)$. By uniqueness of the limit in the weak topology of $H^k_{0,K}(\Omega)$, we must have $h=u$, concluding that finally $u\in  H^k_{0,K}(\Omega)$ as desired.\\
 
 $\bullet$ {Proof of $(4).$}  This item is a direct consequence of (1).
\end{proof}

 
 As seen in the proof of Lemma \ref{QeLemma0}, if $\Omega$ is a Lipschitz domain and $\mathcal{E}$ is an extension operator for $W^{1,2}(\Omega)$, then it holds true that for every Sobolev function $u\in W^{1,2}(\Omega)$,
 
\begin{equation} \label{lim0}
\lim_{r\to 0} \left( \fint_{B_r(x_0) \cap \Omega}| u| \; dx\right) =   \lim_{r\to 0} \left( \fint_{B_r(x_0)} |u| \; dx\right) \quad \text{ for }(1,2)-\text{q.e. }  x_0 \in \partial \Omega. 
\end{equation}
 
 In the paper we would need a more uniform result of the same kind, but where $u$ is not fixed but could also depend on $r$ and converges strongly in $W^{1,2}$. Since we could not find in the literature a short proof for this property, we write in this appendix a complete argument.  For that purpose we use  an explicit extension operator, i.e. a linear mapping $\mathcal{E}:W^{1,2}(\Omega)\to W^{1,2}(\mathbb{R}^2)$ such that $\mathcal{E}(u)=u$ in $\Omega$ and 
\begin{eqnarray}
\|\mathcal{E}(u)\|_{W^{1,2}(\mathbb{R}^2)}\leq C\|u\|_{W^{1,2}(\Omega)}. \label{Extension}
\end{eqnarray}
There are several ways of constructing such an operator. Here we follow  the approach in \cite{EG} as follows. If $\partial \Omega$ locally coincides around $x_0 \in \partial \Omega$ with the graph of the Lipschitz mapping $\gamma:\mathbb{R}^{N-1}\to \mathbb{R}$, and if $u$ is  compactly supported in a cylinder of size $h>0$ centered at $x_0$, one can use the direct formula
\begin{eqnarray}
\mathcal{E}(u)(x',x_N)=
\left\{
\begin{array}{cc}
u(x',x_N) & \text{ if } x_N> \gamma(x') \\
u(x',2\gamma(x')-x_N) & \text{ if } x_N< \gamma(x').
\end{array}
\right. \label{defExxt}
\end{eqnarray}
For the general case  one can cover the boundary $\partial \Omega$ with a finite union of cylinders of same size $h$ and use a partition of unity, leading to an extension operator of the form
\begin{eqnarray}
\mathcal{E}(u)=\sum_i \theta_i \mathcal{E}_i(u), \label{defEG}
\end{eqnarray}
where $\mathcal{E}_i$ is of the form \eqref{defExxt} in a suitable local coordinate system. We refer to Theorem  1 in section 4.4. of \cite{EG} for further detail, where in particular the estimate \eqref{Extension} is established.

Now we focus on the following proposition, and its immediate corollary which is used only in the proof of Lemma \ref{avg-to-0}.

\begin{prop} \label{propositionNice} Let $\Omega$ be a  bounded Lipschitz domain and $\mathcal{E}$ be the extension operator defined in \eqref{defEG}. Let $r\to 0$ and  $(u_r)$ be a sequence in $W^{1,2}(\Omega)$ such that $u_r \to u \in W^{1,2}(\Omega)$. Let   $x_0\in \partial \Omega$. Then 
 $
 \displaystyle{\lim_{r\to 0} \fint_{B(x, r)} \mathcal{E}(u_r) (y) \, dy =0}
 $
 implies 
 $\displaystyle{ \lim_{r\to 0} \fint_{B(x, r) \cap \Omega}  u_r (y) \, dy}=0$. 
 \end{prop}

 \begin{rem}As an immediate corollary, and under the same assumptions of Proposition \ref{propositionNice}, we have 
 $$
 \displaystyle{\lim_{r\to 0} \fint_{B(x, r)}| \mathcal{E}(u_r) (y)| \, dy =0}
 \implies 
 \displaystyle{ \lim_{r\to 0} \fint_{B(x, r) \cap \Omega}  |u_r (y)| \, dy}=0.
 $$ Compare with \eqref{lim0}.
 \end{rem}

\begin{proof} We  first assume that $\partial \Omega \cap B(x_0,r_0)$ coincides with a portion of the graph of  a Lipschitz function $\gamma :\R^N\to \R$, that   $\Omega \cap B(x_0,r)=B(x_0,r)\cap \{ (x',x_N) \; :\; x_N> \gamma(x')\}$ for all $r<r_0$ and that $\mathcal{E}(u)=u(x',2\gamma(x')-x_N)$  in $\{x_N< \gamma(x')\}$. Then we define, for $r<r_0$,
$$B_r^+ = B(x_0,r)\cap \{x_N>\gamma(x')\} \quad \text{ and } \quad  B_r^- = B(x_0,r)\cap \{x_N<\gamma(x')\},$$
and write
$$\int_{B_r(x_0) }\mathcal{E}(u) \; dx  =  \int_{B_r^+ }u \; dx +\int_{B_r^- }u(x' , 2\gamma(x')-x_N) \; dx'dx_N.$$

On the other hand by a simple change of variable wich is linear in the $x_N$ variable, together with Fubini's Theorem, we see that 
$$\int_{B_r^- }|u(x' , 2\gamma(x')-x_N)| \; dx'dx_N = \int_{A_r} |u(x)| \;dx,$$
where $A_r$ is the ``reflected'' domain $\Phi(B_r^-)$, with $\Phi(x',x_N):=(x',2\gamma(x')-x_N)$. Since $\gamma$ is Lipschitz with a constant depending only on $\Omega$, we infer that there exists $\Lambda,\lambda>0$ such that $B_{\lambda r}^+\subset A_r\subset  B_{\Lambda r}^+$ for all $r\leq r_0$. We can therefore estimate the difference between the average in $A_r$ and the one in $B_R^+$ in the following way, denoting by $R=\Lambda r$, (the constant $C$ below depends only on $\Omega$ and can change from line to line)
\begin{eqnarray}
\left|\fint_{A_r} u \;dx -  \fint_{B_{R}^+} u \;dx \right| &=& \left|\fint_{A_r} \left(u  -  \fint_{B_R^+} u\right) \;dx \right|  \leq   \fint_{A_r} \left|u  -  \fint_{B_R^+} u \right| \;dx  \notag \\
&\leq&   Cr^{-2}\int_{B_R^+} \left|u  -  \fint_{B_R^+} u \right| \;dx  \leq  Cr^{-1} \left(\int_{B_R^+} \left(u  -  \fint_{B_R^+} u \right)^2  \;dx  \right)^{\frac{1}{2}} \notag \\
&\leq&  C  \left(\int_{B_R^+}|\nabla u|^2 \;dx  \right)^{\frac{1}{2}}. \label{estimationPoinca}
\end{eqnarray}
For the last inequality we have used the fact that the domains $B_R^+$ are uniformly Lipschitz as $r$ is going to $0$, in order to apply a Poincar\'e inequality with uniform constant in all of the $B_R^+$. Now we apply the above to  $u_r \in W^{1,2}(\Omega)$, instead of $u$ and we know that 
$$\int_{B_{R}^+}|\nabla u_r|^2 \;dx  \xrightarrow{R \to 0} 0,$$
because $u_r\to u$ strongly in $W^{1,2}(\Omega)$, and this proves that 

\begin{eqnarray}\left|\fint_{A_r} u_r \;dx -  \fint_{B_{R}^+} u_r \;dx \right| \xrightarrow{r\to 0} 0.\label{bobo} \end{eqnarray}

Now assume that $\displaystyle{\fint_{B(x_0,r) }\mathcal{E}(u_r)\; dx \xrightarrow{r\to 0} 0}$. Notice that $|A_r|=|B_r^-|$ because the jacobian of $\Phi$ is equal to 1.  Then
\begin{eqnarray}
2\fint_{B_r^+ }u_r \; dx &=&\frac{1}{|B_r^+|}\left(\int_{B_r^+ }u_r \; dx +\int_{A_r }u_r \; dx +\int_{B_r^+ }u_r\; dx- \int_{A_r }u_r \; dx \right) \notag \\
&=&  \frac{|B(x_0,r)|}{|B_r^+|} \; \fint_{B(x_0,r) }\mathcal{E}(u_r) \; dx + \frac{|B_r^-|}{|B_r^+|}(\fint_{B_r^+ }u_r \; dx- \fint_{A_r }u_r \; dx ) + (1-\frac{|B_r^-|}{|B_r^+|})\fint_{B_r^+} u_r \;dx \notag
\end{eqnarray}
So finally
\begin{eqnarray}
(1+\frac{|B_r^-|}{|B_r^+|})\fint_{B_r^+ }u_r \; dx= \frac{|B(x_0,r)|}{|B_r^+|} \; \fint_{B(x_0,r) }\mathcal{E}(u_r) \; dx + \frac{|B_r^-|}{|B_r^+|}(\fint_{B_r^+ }u_r \; dx- \fint_{A_r }u_r\; dx) \xrightarrow{r\to0} 0 \notag
\end{eqnarray}
because of \eqref{bobo}. This  proves that $\fint_{B(x_0,r)\cap \Omega} u_r \; dx =\fint_{B_r^+}u_r \; dx \to 0$ and finishes the proof of the proposition in the particular case when $\mathcal{E}(u)$ coincides with the formula $u(x',2\gamma(x')-x_N)$ under the graph of the Lipschtiz function $\gamma$. \\

In the general case $\partial \Omega$ is covered by a finite number of Lipschitz graphs $\gamma_i$ and $\mathcal{E}(u)$ is of the form $\sum_{i} \theta_i \mathcal{E}_i(u)$ where $\theta_i \in C^\infty_c(\R^2)$ is a partition of unity and $\mathcal{E}_i$ is the extension of $u$ relatively to the graph $\gamma_i$. Applying the above argument to each of the $\mathcal{E}_i(u)$ we conclude as follows.

We   assume that  
$$ \lim_{r\to 0} \left( \fint_{B_r(x_0) }\mathcal{E}(u_r)\; dx\right) =0,$$
which in other words says 
$$ \lim_{r\to 0} \left( \fint_{B_r(x_0) }\sum_{i}\theta_i \mathcal{E}_i(u_r)\; dx\right) =0,$$
or again after change of variable,
$$ \lim_{r\to 0} \fint_{B_r^+} u_r \;dx +\sum_{i}\left( \fint_{A_r^i }\theta_i\circ \Phi_i^{-1} \; u_r\; dx\right) =0,$$
 where $A_r^i$ is the reflexion of $B_r^-$ relatively to the graph of $\gamma_i$, i.e. $A_r^i=\Phi_i(B_r^-)$ with $\Phi_i(x',x_N):=(x',2\gamma_i(x')-x_N)$ (here the coordinate system $(x',x_N)$ should be taken relatively to the one associated with  $\gamma_i$). Now we denote by $I$ the set of indices $i$ for which   $\theta_i$ is not identically zero around $x_0$, and we let $N$ be the cardinal of $I$ (which should actually be at most two).

 Now by construction $\Phi_i(x_0)=x_0$ for all $i$ and $\Phi_i$ is Lipschitz with constant depending only on $\Omega$. Since $\theta_i$ is a smooth function we deduce that 
 \begin{eqnarray}
 |\theta_i\circ \Phi_i^{-1}(x)-\theta_i(x_0)|\leq C|x-x_0|. \label{Lopp}
 \end{eqnarray} 
Let $\Lambda>0$ be such that $A_r^i \subset B_{\Lambda r}^+$ for all small $r$. Then since  $u_r \in L^{2}(\Omega)$ we can  use \eqref{Lopp} to obtain
 \begin{eqnarray}
 \left( \fint_{A_r^i }   |\theta_i\circ \Phi_i^{-1}(x)-\theta_i(x_0)| |u_r| \;dx\right) \leq C  \left(\int_{B(x_0,\Lambda r)} |u_r|^2 \;dx\right)^{\frac{1}{2}} \xrightarrow{r\to 0} 0, \label{Lipp}
 \end{eqnarray}
because $u_r$ converges in $L^2$ to $u$.
Next, following an argument similar to the one of \eqref{estimationPoinca}, denoting by $R:=\Lambda r$, using \eqref{Lipp} and the fact that $\sum_i \theta_i(x_0)=1$ we  obtain some function $e(r)\to 0$  for which
 
\begin{eqnarray}
\left| \fint_{B_r^-} \mathcal{E}( u_r) \;dx -  \fint_{B_{R}^+} u_r \;dx \right|  &=&\left|\left(\sum_{i\in I}\fint_{A_r^i} \theta_i\circ \Phi_i^{-1} \; u_r \;dx\right) -  \fint_{B_{R}^+} u_r \;dx \right|  \notag \\
&\leq & \left|\sum_i \theta_i(x_0)\fint_{A_r^i} \left(u_r  -  \fint_{B_R^+} u_r\right) \;dx \right| +e(r)  \notag \\
&\leq &  Cr^{-2}\int_{B_R^+} \left|u_r  -  \fint_{B_R^+} u_r \right| \;dx + e(r) \notag \\
&\leq & C  \left(\int_{B_R^+}|\nabla u_r|^2 \;dx  \right)^{\frac{1}{2}} + e(r) \notag\\
& & \xrightarrow{r\to 0}0. \label{estimationPoinca2}
\end{eqnarray}

    Then as before we can write
    
\begin{eqnarray}
(1+\frac{|B_r^-|}{|B_r^+|})\fint_{B_r^+ }u_r \; dx= \frac{|B(x_0,r)|}{|B_r^+|} \; \fint_{B(x_0,r) }\mathcal{E}(u_r) \; dx + \frac{|B_r^-|}{|B_r^+|}(\fint_{B_r^+ }u_r \; dx- \fint_{B_r^- }\mathcal{E}(u_r) \; dx) \xrightarrow{r\to0} 0 \notag
\end{eqnarray}
which finishes the proof in the general case, and thus the proposition is now proved.
\end{proof}

 }

 We end this section with  variants of the Poincar\'e inequality related to the traces of functions on positive capacity subsets. The first result can be compared to \cite[Corollary 8.2.2]{ah}.
 
 \begin{cor}\label{qe-poinc}
Let $\Omega \subset \R^2$ be a bounded Lipschitz domain and let  the closed set $K\subset \overline\Omega$ be such that ${\rm Cap}_{1,2} (K)>0$. Then there exists a constant $C>0$ such that for all {$u\in H^1_{0,K}(\Omega)$,}
$$
\|u\|_{L^2(\Omega)} \le C \|\nabla u\|_{L^2(\Omega)},
$$
{where $H^k_{0,K}(\Omega)$ is the space defined in Lemma \ref{QeLemma0}.}
 \end{cor}
 
 \begin{proof}
 
We prove the result by contradiction. Assume that there exists a sequence $u_k \in H^1_{0,K}(\Omega)$ such that 
$$
\|u_k\|_{L^2(\Omega)} \geq k \|\nabla u_k\|_{L^2(\Omega)}.
$$   By renormalizing the sequence, we can assume that $\|\nabla u_k\|_{L^2(\Omega)} \to 0$, as $k\to \infty$, while $\|u_k\|_{L^2(\Omega)} =1$.  It follows that, passing to a subsequence, $u_k$ converges strongly to a constant $a_0\neq 0$.  But  it immediately follows from Lemma \ref{QeLemma0}  that $ a_0 \equiv 0$ on $K$ quasi-everywhere, which contradicts $a_0\neq 0$ as a constant.
  \end{proof}
 
 The second result concerns uniform constants for the Poinca\'e inequality on thin domains $\Omega^h$, where the trace of the functions vanish on subsets of the form $K_h \times \{0\} \subset \Omega \times \{0\}$. 
 It should be compared with  \cite[Theorem D.1]{LM11}.
  \begin{theorem}\label{thin-poinc}
  Let $\Omega$ be a bounded Lipschitz domain, $\Omega^h:= \Omega \times (0,1)$, and let 
  $$
  {\bf A}^h_K := \{u\in W^{1,2}(\Omega^h); \,\, u|_{K_h \times \{0\}} =0 \}.
  $$
  Then there exists $h_0>0$ such that for a constant $C>0$ uniform in $h<h_0$  we have 
  $$
 \forall  u \in   {\bf A}^h_K \quad \|u\|_{L^2(\Omega^h)} \le C \|\nabla u\|_{L^2(\Omega^h)}.
  $$
  \end{theorem}
   
We will provide a sketch of the proof,  based on \cite[Theorem D.1]{LM11}, adapted to our situation.  We will follow the steps leading to Lemma \ref{mollified-grady}, replacing the rigidity estimate in Corollary \ref{fixed-Id} by the Poincar\'e inequality applied to the Lipschitz domains:
$$
\|u - a_{x'}\|_{L^2({\mathcal Q}_j (x', h))} \le Ch^2 \|\nabla u\|_{L^2({\mathcal Q}_j (x', h))}
$$
with $a_{x'}=0$ whenever $K_{ch} \cap \overline{\Phi_j^{-1}(S_{\xi', h})} \neq \emptyset$. We can then conclude with the following lemma parallel to Lemma \ref{tildF'}:
  \begin{lemma}\label{a-poinc}
 Let $\Omega$, $\Omega^h$, $K$ be as defined above. Then there exist constants $h_0>0$, $C>0$, $0<\bar c<1$, depending only on $\Omega$ and $K$, such that given $h<h_0$, $u\in \mathbf{A}^h_K$, there exists a   scalar function $a:\Omega \to \R$  which  vanishes on $K_{\bar ch}$ such that the estimates  
    \begin{equation}\label{for-a-poinc}
  \ds   \|u- a\|^2_{L^2(\Omega^h)}  \le C h^2  \|\nabla u\|^2_{L^2(\Omega^h)}, \quad   \|\nabla a\|^2_{L^2(\Omega)} 
    \le \frac{C}{h} \|\nabla u\|^2_{L^2(\Omega^h)}, 
  \end{equation} and 
  \begin{equation}\label{a-infty}
  \|a\|^2_{L^\infty(\Omega)} \le \frac C{h}  \|\nabla u\|^2_{L^2(\Omega^h)}, 
\end{equation}  hold true.
 \end{lemma} 
 
\noindent {\bf Proof of Theorem \ref{thin-poinc}.} Let $u \in  \mathbf{A}^h_K$ and let $a\in W^{1,2}(\Omega)$ be chosen according to Lemma \ref{a-poinc}. Note that $a\in H^1_{0,K}(\Omega)$, and hence by Corollary \ref{qe-poinc} we obtain for a uniform constant $C$
 $$
 \|a\|^2_{L^2(\Omega)} \le C\|\nabla a\|^2_{L^2(\Omega)},
 $$ which yields, through  the second inequality in \eqref{for-a-poinc}
 $$
  \|a\|^2_{L^2(\Omega^h)}= h\|a\|^2_{L^2(\Omega)} \le Ch \|\nabla a\|^2_{L^2(\Omega)} \le C\|\nabla u\|^2_{L^2(\Omega^h)} . 
 $$ The conclusion follows from the first inequality in \eqref{for-a-poinc}.
 
 
 
\bibliography{bib1}
\bibliographystyle{plain}

 





\end{document}